\documentstyle[amssymb,amstex,a4wide,12pt]{article}
\newtheorem{theorem}{Theorem}
\newtheorem{lemma}{Lemma}
\newtheorem{corollary}{Corollary}

\newcommand{\Edges}{\operatorname{Edges}}
\newcommand{\Vertices}{\operatorname{Vertices}}
\newcommand{\Angles}{\operatorname{Angles}}
\newcommand{\iin}{\operatorname{in}}
\newcommand{\oout}{\operatorname{out}}
\newcommand{\card}{\operatorname{card}}
\newcommand{\conv}{\operatorname{conv}}
\newcommand{\lnaw}{\langle\langle}
\newcommand{\pnaw}{\rangle\rangle}
\newcommand{\lgth}{\ell}
\newcommand{\lodc}{\overline{[[}}
\newcommand{\podc}{\overline{]]}}
\newcommand{\e}{\operatorname{e}}
\newcommand{\An}{\operatorname{An}}
\newcommand{\innt}{\operatorname{Int}}
\newcommand{\diam}{\operatorname{diam}}
\newcommand{\dist}{\operatorname{dist}}

\begin{document}

\setlength{\baselineskip}{1.3\baselineskip}

\author{Tomasz Schreiber\footnote{Mailing address: Tomasz Schreiber, Faculty of
 Mathematics \& Computer Science, Nicolaus Copernicus University,
 ul. Chopina 12 $\slash$ 18, 87-100 Toru\'n, Poland; tel.: (++48) (+56)
 6112951, fax: (++48) (+56) 6228979; e-mail: {\tt tomeks at mat.uni.torun.pl}},\\
        Faculty of Mathematics \& Computer Science,\\
        Nicolaus Copernicus University, Toru\'n, Poland.}

\title{Non-homogeneous polygonal Markov fields in the plane: graphical
       representations and geometry of higher order correlations}
\date{}
\maketitle

\paragraph{Abstract} {\it We consider polygonal Markov fields originally introduced by
 Arak and Surgailis \cite{AS1}. Our attention is focused on fields with nodes of 
 order two, which can be regarded as continuum ensembles of non-intersecting
 contours in the plane, sharing a number of features with the two-dimensional
 Ising model. We introduce non-homogeneous version of polygonal fields in
 anisotropic enviroment. For these fields we provide a class of new graphical
 constructions and random dynamics. These include a generalised dynamic
 representation, generalised and defective disagreement loop dynamics as well
 as a generalised contour birth and death dynamics. Next, we use these
 constructions as tools to obtain new exact results on the geometry of
 higher order correlations of polygonal Markov fields in their consistent
 regime.} 


\section{Introduction}
 The first example of a polygonal Markov field has been provided by Arak \cite{A1}.
 Further developments are due to Arak \& Surgailis \cite{AS1,AS2}, Surgailis
 \cite{SU1} and Arak, Clifford \& Surgailis \cite{ACS}. The fields with nodes of order $2,$ 
 or V-shaped nodes for short, which are in the focus of our attention in this paper,
 arise as ensembles of self-avoiding closed polygonal contours in the plane
 interacting by hard core exclusions, possibly with some further terms entering
 the Hamiltonian. Under a particular choice of the Hamiltonian the polygonal fields
 enjoy striking properties including consistency (the field constructed in a
 subdomain $D \subseteq D' \subseteq {\Bbb R}^2$ coincides with the restriction
 to $D$ of the field constructed in $D'$) as well as availability of an explicit formula
 for the partition function and some other numerical characteristics of the field,
 see \cite{AS1}. Under even milder conditions one can guarantee
 the isometry invariance of the field as well as the two-dimensional germ-Markov
 property stating that the conditional behaviour of the field in an open bounded
 domain depends on the exterior configuration only through arbitrarily close
 neighbourhoods of the boundary, see ibidem. A particularly interesting class
 of processes seem to be length- and area-interacting modifications of consistent
 fields, which not unexpectedly exhibit many features analogous to those of the
 two-dimensional Ising model, including the presence of first order phase 
 transition at low enough temperatures (Nicholls, \cite{N1}; Schreiber, \cite{SC05})
 and low-temperature phase separation phenomenon (Schreiber, \cite{SC06}). Consistent
 polygonal fields and their length-interacting modifications have also interesting
 statistical applications where they are used as priors in Bayesian image analysis
 (Clifford \& Nicholls, \cite{CN94}; Kluszczy\'nski, van Lieshout \& Schreiber,
 \cite{KLS05,KLS07}; van Lieshout \& Schreiber, \cite{LS07,SL08}).

 With all the afore-mentioned similarities the polygonal fields enjoy with the
 Ising model as well as with more general lattice models used in Bayesian image analysis,
 one salient difference is that whereas the latter can be analysed using a
 rich mixture of combinatorial, analytic and geometric techniques,
 the overwhelming majority of crucial results for the polygonal Markov fields
 have been obtained via purely geometric methods usually going under the
 guise of graphical constructions. These were first provided by Arak 
 \& Surgailis \cite{AS1,AS2} and Arak, Clifford \& Surgailis \cite{ACS} in
 the form of the so-called {\it dynamic representations} for consistent
 polygonal fields. Later, we have introduced {\it disagreement loop dynamics} 
 and {\it contour birth and death representation} for broader classes of
 polygonal fields, see \cite{SC05}. The purpose of the present paper is twofold.
 \begin{itemize}
  \item First, we introduce a class of new graphical constructions and random
        dynamics for a rather general class of non-homogeneous polygonal Markov
        fields with V-shaped nodes. These constructions include a generalised dynamic
        representation, generalised and defective disagreement loop dynamics as well
        as a generalised contour birth and death dynamics.
  \item Next, we use these constructions as tools to obtain new exact results
        on the geometry of higher order correlations of polygonal Markov fields
        in the consistent regime.
  \end{itemize}
  We envision to use the new graphical representations in our current research
  in progress on polygonal fields, especially in context of obtaining further
  exact formulae for higher order correlations and tentatively also to obtain
  exact formulae for the free energy of the field (exact solution) at least
  in the rectangular case. Due to their intrinsically algorithmic nature 
  the new constructions will also be used for Bayesian image segmentation
  purposes, which is our work in progress in cooperation with van Lieshout,
  see Schreiber \& van Lieshout \cite{SL08} for some partial developments in
  the particular tesselation-based set-up.  
 
  The paper is organised as follows. In Section \ref{NHPMF} below we introduce 
  general non-homogeneous polygonal fields in non-homogeneous and anisotropic
  environments. This is done by imposing a general {\it activity measure} on
  the space of straight lines in the plane, determining how likely the presence
  of an edge along a given line is. The usual Ising-type length interaction
  is replaced by its anisotropic version there \--- the energy cost of
  producing an edge crossing lines of high activities is higher than the
  cost of an edge of the same length crossing only low activity lines.
  Further, in Section \ref{CONSIS} general properties of such fields in
  their {\it consistent} regime are discussed and a non-homogeneous counterpart
  of the Arak \& Surgailis \cite{AS1} dynamic representation is developed. 
  Next, in Section \ref{GDR} we develop a generalised dynamic representation
  for consistent polygonal fields, which allows for a flexible choice of the field
  creation dynamics, parametrised by an increasing family of convex compact
  subsets of the field domain. In the subsequent Section \ref{CORCR} we employ
  the generalised dynamic representation to provide exact factorisation statements
  for higher order correlation functions of polygonal fields in consistent regime.
  In the next Section \ref{DLOOP} we introduce variants of the so-called disagreement
  loop dynamics for polygonal fields, including the generalised and defective
  disagreement loop dynamics. These are then used in Section \ref{CORRDEC}
  to establish exponential decay of higher order correlations in the particular
  set-up of rectangular fields. Finally, in Section \ref{CBDD} we introduce
  general contour birth and death dynamics for polygonal fields in low
  enough temperature, in the spirit of Fern\'andez, Ferrari \& Garcia
  \cite{FFG1,FFG2,FFG3} and yielding a perfect simulation scheme directly
  from the thermodynamic limit.    


\section{Non-homogeneous polygonal Markov fields}\label{NHPMF}
 For an open bounded convex set $D$ define the family $\Gamma_D$ of admissible polygonal
 configurations on $D$ by taking all the finite planar graphs $\gamma$ in $D \cup \partial D,$
 with straight-line segments as edges, such that
 \begin{description}
  \item{\bf (P1)} the edges of $\gamma$ do not intersect,
  \item{\bf (P2)} all the interior vertices of $\gamma$ (lying in $D$) are of degree $2,$
  \item{\bf (P3)} all the boundary vertices of $\gamma$ (lying in $\partial D$) are of degree $1,$
  \item{\bf (P4)} no two edges of $\gamma$ are colinear.
 \end{description}
 In other words, $\gamma$ consists of a finite number of disjoint polygons, possibly nested
 and chopped off by the boundary. Further, for a finite collection $(l) = (l_i)_{i=1}^n$ of
 straight lines intersecting $D$ we write $\Gamma_D(l)$ to denote the family of admissible
 configurations $\gamma$ with the additional properties that
 $\gamma \subseteq \bigcup_{i=1}^n l_i$ and $\gamma \cap l_i$ is a single interval
 of a strictly positive length for each $l_i, i=1,...,n,$ possibly with some
 isolated points added.

 We shall also consider the subfamily $\Gamma_{D|\emptyset} \subseteq \Gamma_D$
 of empty boundary polygonal configurations in $D$ consisting of those $\gamma \in \Gamma_D$  
 which have no boundary vertices. For $(l)$ as above we put $\Gamma_{D|\emptyset}(l) 
 := \Gamma_D(l) \cap \Gamma_{D|\emptyset}.$ 

 For a Borel subset of $A \subseteq {\Bbb R}^2$ by $[[A]]$ we shall denote the
 family of all straight lines hitting $A$ so that in particular $[[{\Bbb R}^2]]$
 stands for the collection of all straight lines in ${\Bbb R}^2.$ We shall also
 write $\lodc A \podc$ for the family of all linear segments in
 ${\Bbb R}^2$ hitting $A.$   
 Consider a non-negative Borel measure ${\cal M}$ on the $[[{\Bbb R}^2]]$
 such that
 \begin{description}
  \item{\bf (M1)} ${\cal M}([[A]]) < \infty$ for all bounded Borel $A \subseteq {\Bbb R}^2,$
  \item{\bf (M2)} ${\cal M}([[\{x\}]]) = 0$ for all $x \in {\Bbb R}^2.$
 \end{description}
 Below, the measure ${\cal M}$ will be interpreted as the activity measure on
 $[[{\Bbb R}^2]].$ Let $\Lambda^{\cal M}$ be the Poisson line process on $[[{\Bbb R}^2]]$
 with intensity measure ${\cal M}$ and write $\Lambda^{\cal M}_D$ for its restriction to
 the domain $D.$ Further, define the Hamiltonian 
 $L^{\cal M} : \Gamma_D \to {\Bbb R}_+$ given by 
 \begin{equation}\label{HAMILT}
  L^{\cal M}(\gamma) := \sum_{e \in \Edges(\gamma)} {\cal M}([[e]]),\;\; \gamma \in \Gamma_D.
 \end{equation}
 We argue that the energy function $L^{\cal M}$ should be regarded as
 an anisotropic  environment-specific
 version of the length functional. Indeed, interpreting the activity ${\cal M}(dl)$ of
 a line $l$ hitting an edge $e \in \Edges(\gamma)$ as the likelihood of a new edge being
 created along $l$ intersecting and hence fracturing the edge $e$ in $\gamma,$
 we note that, roughly speaking, the value of ${\cal M}([[e]])$ determines how likely the
 edge $e$ is to be fractured by another edge present in the environment.
 In other words, $L^{\cal M}(\gamma)$ determines {\it how difficult it is to
 create} the whole graph $\gamma \in \Gamma_D$ without fractures in the environment
 ${\cal M}$ \--- note that due to the anisotropy of the environment there may be graphs
 of a higher (lower) total edge length than $\gamma$ and yet of lower (higher) energy
 and thus easier (more difficult) to create and to keep unfractured due to the lack
 (presence) of high activity lines likely to fracture their edges.
   
 With the above notation, for $\beta \in {\Bbb R}$ further referred to as the
 inverse temperature (for mathematical convenience we admit also the unphysical
 negative values of inverse temperatures),
 we define the polygonal field ${\cal A}^{{\cal M};\beta}_D$ in $D$ with activity
 measure ${\cal M}$ to be the Gibbsian modification of the process induced on
 $\Gamma_D$ by $\Lambda^{\cal M}_D,$ with the Hamiltonian $L^{\cal M}$ at 
 inverse temperature $\beta,$ that is to say
 \begin{equation}\label{GIBBS1}
  {\Bbb P}\left( {\cal A}^{{\cal M};\beta}_D \in G \right)
  := \frac{{\Bbb E} \sum_{\gamma \in \Gamma_D(\Lambda^{\cal M}_D) \cap G}
     \exp\left(- \beta L^{\cal M}(\gamma) \right)}
          {{\Bbb E} \sum_{\gamma \in \Gamma_D(\Lambda^{\cal M}_D)}
     \exp\left(- \beta L^{\cal M}(\gamma) \right)}
 \end{equation}
 for all sets $G \subseteq \Gamma_D$ Borel measurable with respect to, say, the usual
 Hausdorff distance topology. Note that this definition can be rewritten as
 \begin{equation}\label{GIBBS2}
    {\Bbb P}({\cal A}^{{\cal M};\beta}_D \in d\gamma) \propto 
    \exp(- \beta L^{\cal M}(\gamma)) \prod_{e \in \Edges(\gamma)}
    {\cal M}(dl[e]),\; \gamma \in \Gamma_D,
 \end{equation}
 where $l[e]$ is the straight line extending $e.$ In other words, the probability
 of having ${\cal A}^{{\cal M};\beta}_D \in d\gamma$ is proportional to 
 the Boltzmann factor $\exp(-\beta L^{\cal M}(\gamma))$
 times the product of edge activities ${\cal M}(dl[e]),\; e \in \Edges(\gamma).$
 Observe also that this construction should be regarded as a specific version
 of the general polygonal model given in (2.11) in \cite{AS1}.
 The finiteness of the partition function
 \begin{equation}\label{PAFA}
   {\cal Z}^{{\cal M};\beta}_D := {\Bbb E} \sum_{\gamma \in \Gamma_D(\Lambda^{\cal M}_D)}
      \exp\left(- \beta L^{\cal M}(\gamma) \right) < \infty
 \end{equation}
 for all $\beta \in {\Bbb R}$ is not difficult to verify, see (\ref{ZZ2}) below.
 We shall also consider the empty boundary version of the above
 construction,  with ${\cal A}^{{\cal M};\beta}_{D|\emptyset}$ arising in law as 
 ${\cal A}^{{\cal M};\beta}_D$ conditioned on staying within $\Gamma_{D|\emptyset},$
 that is to say not containing boundary vertices. It is easily seen that alternatively
 ${\cal A}^{{\cal M};\beta}_{D|\emptyset}$ can
 be defined by rewriting (\ref{GIBBS1}) above with $\Gamma_D$ replaced with
 $\Gamma_{D|\emptyset}$ throughout. The corresponding partition function
 ${\cal Z}^{{\cal M};\beta}_{D|\emptyset}$ is always finite as well,
 see (\ref{ZZ1}) below. 

 In the sequel we also consider the thermodynamic limits for such fields,
 which take their values in the space $\Gamma := \Gamma_{{\Bbb R}^2}$ of
 whole-plane admissible configurations, with obvious meaning of this
 notation. Note in
 this context that $\Gamma_{{\Bbb R}^2} = \Gamma_{{\Bbb R}^2|\emptyset}.$  

 As a particular case of this general construction we shall consider 
 {\it rectangular fields}. Slightly abusing the terminology,
 by a rectangular field we shall mean a polygonal field 
 given by (\ref{GIBBS1}) and (\ref{GIBBS2}) whose activity measure
 ${\cal M}$ concentrates on translates of two non-parallel straight
 lines \--- in other words, its edges can follow precisely two different
 directions. Clearly, upon a non-singular affine transform such a
 rectangular field can be  made into a process where all angles between
 edges are straight, whence the name. Assuming with no loss of generality
 that the two  admissible directions for the edges of the field are
 parallel to the coordinate axes, we can further transform
 each rectangular field into the {\it standard rectangular field} where the 
 measure ${\cal M}$ is given by ${\cal M}(dl)$ being $dh$ if $l$ is a 
 parallel translate of a coordinate axe by distance $h,\; h \in
 {\Bbb R},$ and $0$ otherwise. Indeed, putting
 $$ F_{\leftrightarrow}(x) := \left\{ \begin{array}{cl} {\cal M}((0,x] \times \{ 0 \}), & \mbox{ if }
    x \geq 0,\\
    - {\cal M}((x,0] \times \{ 0 \}), & \mbox{ otherwise } \end{array} \right. $$
 and 
 $$ F_{\updownarrow}(y) := \left\{ \begin{array}{cl} {\cal M}(\{ 0 \} \times (0,y]), & \mbox{ if }
    y \geq 0,\\
    - {\cal M}(\{ 0 \} \times (y,0]), & \mbox{ otherwise, } \end{array} \right. $$
 it is readily verified that the map  ${\Bbb R}^2 \ni (x,y) \mapsto
 (F_{\leftrightarrow}(x),F_{\updownarrow}(y))$ does define such a transformation.
 Note that both $F_{\leftrightarrow}$ and $F_{\updownarrow}$ are finite by {\bf (M1)}
 and continuous by {\bf (M2)}.
 Although in general this transformation is not a bijection, it is easily seen
 to be bijective on the support of ${\cal A}^{\cal M},$ with the inverse defined by
 $(s,t) \mapsto (\sup \{ x,\; F_{\leftrightarrow}(x) \leq s \},
  \sup \{ y,\; F_{\updownarrow}(y) \leq t\}).$

 Even though in general we admit all $\beta \in {\Bbb R},$ in this paper
 we shall mainly study two principal temperature regions. The first
 one is the {\it consistent regime} corresponding to $\beta = 1$ and 
 introduced in  Section \ref{CONSIS} \---  we shall argue
 that this particular choice of temperature parameter places us in the
 context of a non-homogeneous version of Arak-Surgailis \cite{AS1} construction
 for the so-called {\it consistent} polygonal fields, thus ensuring the
 availability of an appropriate dynamic representation for our process in
 terms of equilibrium evolution of one-dimensional particle systems tracing
 the boundaries of the field in two-dimensional space-time, as discussed in detail
 in Subsection \ref{DYNRES} below, see also \cite{AS1}. The afore-mentioned consistency
 property, arising as a consequence of the dynamic representation and further
 discussed in Subsection \ref{CDYNRES}, means here that for each open bounded
 convex $D \subseteq D' \subset {\Bbb R}^2$ the field
 ${\cal A}^{{\cal M};1} := {\cal A}^{{\cal M};1}_D$ coincides in law with
 ${\cal A}^{{\cal M}}_{D'} := {\cal A}^{{\cal M};1}_{D'}
 \cap D,$ thus allowing for a direct construction of the infinite volume process
 (thermodynamic limit) ${\cal A}^{\cal M} := {\cal A}^{{\cal M};1}$ on the whole
 ${\Bbb R}^2.$ A number of further properties can be concluded from the dynamic
 representation for $\beta = 1,$ including the explicit knowledge of the partition
 function $Z^{{\cal M};1}_D$ as made precise in Theorem \ref{WlasnosciZgd} below.
 The consistent fields are the main object of study in Sections \ref{CONSIS},
 \ref{GDR}, \ref{CORCR} and \ref{CORRDEC}, and play a crucial role also in
 Section \ref{DLOOP}.

 The second regime in the focus of our interest is the low temperature region
 (large positive $\beta$) where long range point-to-point correlations are present,
 giving rise to the spontaneous magnetisation phenomenon, see Corollary
 \ref{FinContNest} in Section \ref{CBDD} and the discussion below it.
 All disagreement loop dynamics developed in Section \ref{DLOOP} of this paper 
 admit their versions working for the low temperature region.
 Moreover, in Section \ref{CBDD} we provide another important tool characterising
 the behaviour of our model in this regime, namely a {\it contour birth and death
 graphical construction} in the spirit of Fern\'andez, Ferrari and Garcia
 \cite{FFG1,FFG2,FFG3} put together with our random walk representation for the
 so-called {\it free contour measure} specific to the model. Applying this
 graphical construction we are able to establish the existence of the thermodynamic
 limit for ${\cal A}^{{\cal M};\beta}$ as well as certain mixing-type results for
 general activity measures satisfying (\ref{HOMOGENBD}) in analogy with
 Section 4 in \cite{SC05}, see Corollary \ref{BetaMix}.   
 
 As a by-product of our considerations both in the consistent and low temperature regimes
 we obtain efficient exact simulation algorithms for ${\cal A}^{{\cal M};\beta},$ both
 in finite windows and directly from the thermodynamic limit.       


\section{The consistent regime}\label{CONSIS}
 As mentioned above, we show here that the choice $\beta=1$  puts us in the
 so-called {\it consistent regime} of the polygonal field, where it enjoys
 remarkable features arising as consequences of a (non-homogeneous) version
 of the dynamic representation developed by Arak \& Surgailis \cite{AS1},
 see Section 4 there. In Subsection \ref{DYNRES} below, we give the
 non-homogeneous dynamic representation and further in Subsection
 \ref{CDYNRES} we discuss its consequences.   

\subsection{Consistent fields and their dynamic representation}\label{DYNRES}
 To describe the dynamic representation, we interpret the open convex
 domain $D$ as a set of {\it space-time} points $(t,y) \in D,$ with $t$ referred to as
 the {\it time} coordinate and with $y$ standing for the {\it spatial} coordinate of
 a particle at the time $t.$ In this language, a straight line segment in $D$
 stands for a piece of the space-time trajectory of a freely moving particle.
 In the sequel by a particle we shall understand a function $t \mapsto y_t$ 
 assigning to each time moment $t$ the spatial location $y_t$ of the particle
 at the time $t.$ Moreover, by the trajectory of a particle  we shall mean
 the graph of this function, i.e. the set of points $(t,y_t).$ The particle
 trajectories considered in this paper are piecewise linear and hence almost
 everywhere differentiable. Clearly, in the language being introduced,
 the derivative of $y_t$ at $t$ corresponds to the current velocity of the
 particle at the time $t.$ Below, for notational convenience we shall omit
 the subscript $t$ writing $y$ instead of $y_t.$
 For a straight line $l$ non-parallel to the spatial axis and crossing the
 domain $D$ we define in the obvious way its entry point to $D,\;\iin(l,D)
 \in \partial D$ and its exit point $\oout(l,D) \in \partial D.$

 On ${\Bbb R}^2$ we construct the {\it birth measure} $\lnaw{\cal M}\pnaw$
 given for each Borel $A \subseteq {\Bbb R}^2$ by 
 \begin{equation}\label{MM}
  \lnaw {\cal M} \pnaw(A) := 
  \frac{1}{2} {\cal M} \times {\cal M}( \{ (l_1,l_2),\; l_1 \cap l_2 \subset A \}) = 
  {\Bbb E}\card\{ \{l_1,l_2\} \subseteq \Lambda^{\cal M},\;l_1 \neq l_2,\; l_1 \cap l_2 \subset A \}. 
 \end{equation}
 Note in this context that ${\cal M} \times {\cal M}(\{(l,l),\; l \in [[{\Bbb R}^2]]\}) = 0$
 by {\bf (M2)}. Likewise, on $\partial D$ we construct the
 boundary birth measure $\lnaw {\cal M};\partial D \pnaw$ given
 for each Borel $B \subseteq \partial D$ by 
 \begin{equation}\label{MD}
  \lnaw {\cal M};\partial D \pnaw(B) := 
  {\cal M}(\{ l,\; \iin (l,D) \in B \}) = 
  {\Bbb E}\card\{ l \in \Lambda^{\cal M},\; \iin(l,D) \in B \}.
 \end{equation}

 We choose the space-time birth coordinates
 for the new particles according to a Poisson point process in $D$ with intensity
 measure $\lnaw {\cal M} \pnaw$ (interior birth
 sites) superposed with a Poisson point process on the boundary (boundary
 birth sites) with the intensity measure
 $\lnaw {\cal M};\partial D\pnaw.$
 Each interior birth site $x \in D$ emits two particles, moving with initial
 velocities $v'$ and $v''$ chosen according to the joint distribution
 \begin{equation}\label{THETAWPRO}
   \theta^{\cal M}_{x}(dv',dv'') := \frac{
   {\cal M} \times {\cal M}(\{ (l_1,l_2),\; l_1 \cap l_2 \in dx, v[l_1] \in dv', v[l_2] \in dv'' \})}
   {2\lnaw {\cal M} \pnaw(dx)}   
 \end{equation}
 where $v[l]$ is the velocity of a particle whose space-time trajectory is
 represented by the straight line $l.$ Note that this is equivalent to
 choosing the directions of the straight lines representing the space-time
 trajectories of the emitted particles according to the distribution of the
 {\it typical angle} between two lines of $\Lambda^{\cal M}$ at $x,$ that is
 to say the angle arising by conditioning $\Lambda^{\cal M}$ on containing
 two lines intersecting at $dx,$ see also Sections 3 and 4 in \cite{AS1}
 and the references therein. In other words, (\ref{THETAWPRO}) implies
 that the interior birth at $l_1 \cap l_2$ in $D$ of two particles moving along straight
 lines $l_1$ and $l_2$ respectively happens with probability
 ${\cal M}(dl_1) {\cal M}(dl_2).$
 Each boundary birth site $x \in \partial D$
 yields one particle with initial speed $v$ determined according to the distribution
 $\theta^{{\cal M};\partial}_x(dv)$
 identified by requiring that the direction of the line entering $D$ at $x$ and
 representing the space-time trajectory of the emitted particle be chosen according
 to the distribution of a straight line $l \in \Lambda^{\cal M}$ conditioned
 on the event $\{ x = \iin(l,D) \},$ that is to say
 \begin{equation}\label{THETAWPRO2}
  \theta^{{\cal M};\partial}_x(dv) := \frac{ {\cal M}(\{ l,\; \iin(l,D) \in dx, v[l] \in dv\})}
   {\lnaw{\cal M};\partial D\pnaw(dx)}.
 \end{equation}
 Again, this means that the boundary birth at $l \cap \partial D$ of a particle moving along
 a straight line $l$ occurs with probability ${\cal M}(dl).$

 All the  particles evolve independently in time according to the following rules.
 \begin{description}
  \item{\bf (E1)} Between the critical moments listed below each particle
                  moves freely with constant velocity so that $dy = v dt,$
  \item{\bf (E2)} When a particle touches the boundary $\partial D,$ it dies,
  \item{\bf (E3)} In case of a collision of two particles (equal spatial coordinates $y$
        at some moment $t$ with $(t,y) \in D$), both of them die,
  \item{\bf (E4)} The time evolution of the velocity
        $v_t$ of an individual particle is given by the following pure-jump inhomogeneous
        Markov process: denoting by $l_t$ the straight line extending the present segment
        of the space-time trajectory of the particle at time $t$ so that $v_t = v[l_t],$
        we have  
        $$ {\Bbb P}(l_{t+dt} \in dl' \;|\; l_t = l) = 
           {\bf 1}_{\{l'\in [[\,l^{[t,t+dt]}\,]] \}}
           \lnaw{\cal M}\pnaw(dl') 
        $$
        with $l^{[t,t+dt]}$ standing for the segment along $l$ between time
        moments $t$ and $t+dt.$ 
 \end{description}
 The evolution rule {\bf (E4)} can be interpreted as follows: in the course of
 its movement along a short line segment $l^{[t,t+dt]}$ a particle turns and
 starts moving along another straight line $l' \in [[\, l^{[t,t+dt]}\, ]]$
 with probability $\lnaw {\cal M}\pnaw(dl')$ which
 agrees with the interpretation of ${\cal M}$ as the line activity measure as
 discussed above.  

 The following theorem shows that the polygonal field obtained in the course of
 the above dynamic construction is in fact ${\cal A}^{\cal M}_D = 
 {\cal A}^{{\cal M};1}_D,$ in analogy with the results of Arak \& Surgailis \cite{AS1},
 Section 4, for homogeneous fields.
 \begin{theorem}\label{ReprDyn}
  The polygonal field traced by the particle system
  constructed above coincides in distribution with ${\cal A}^{\cal M}_D =
  {\cal A}^{{\cal M};1}_D.$ Moreover, we have
  \begin{equation}\label{ZZ}
   {\cal Z}^{{\cal M};1}_D = \exp(\lnaw {\cal M} \pnaw(D)).
  \end{equation}
 \end{theorem} 

 \paragraph{Proof}
  We pick some $\gamma \in \Gamma_D$ and calculate the
  probability that the outcome of the above dynamic construction falls into $d\gamma.$
  To this end, we note that:
  \begin{itemize}
   \item Each edge $e \in \Edges(\gamma)$ with initial vertex (lower time coordinate) lying
         on $\partial D$ contributes to the considered probability the factor 
         ${\cal M}(dl[e])$ (boundary birth of a particle tracing the edge)
         times $\exp(-{\cal M}([[e]]))$ (no velocity updates along $e$).
   \item Each of the two edges $e_1,e_2 \in \Edges(\gamma)$ steming from a common interior
         birth site $l[e_1] \cap l[e_2]$ yields the factor ${\cal M}(dl[e_i]),\;
         i=1,2,$ (coming from the birth probability) times $\exp(-{\cal M}([[e_i]]))$  
         (no velocity updates along $e_i$),
  \item  Each of the edges $e \in \Edges(\gamma)$ arising due to a velocity update of a particle
         yields the factor ${\cal M}(dl[e])$ (velocity update probability) times
         $\exp(-{\cal M}(dl[e]))$ (no velocity updates along $e$),
  \item  The absence of interior birth sites in $D \setminus \gamma$ yields the factor
         $\exp(-\lnaw {\cal M} \pnaw(D)),$
  \item  Finally, the absence of boundary birth sites at $\partial D \setminus \gamma$ yields
         the additional factor $\exp(-\lnaw {\cal M};\partial D\pnaw(D))
         $ $=$ $\exp(-{\cal M}(D)).$
 \end{itemize}      
 Putting these observations together we conclude that the probability element of $\gamma$
 being traced by the particle system is
 $$ \frac{\exp\left(-L^{\cal M}(\gamma)\right)}
         {\exp(\lnaw{\cal M}\pnaw(D))}
    \frac{\prod_{e\in\Edges(\gamma)} {\cal M}(dl[e])}
         {\exp({\cal M}(D))} = 
    \frac{\exp\left(-L^{\cal M}(\gamma)\right)}
         {\exp(\lnaw{\cal M}\pnaw(D))}        
         d{\Bbb P}(\gamma \in \Gamma_D(\Lambda^{\cal M})), $$
 where the event $\{ \gamma \in \Gamma_D(\Lambda)^{\cal M} \}$ is easily seen
 to correspond to the situation where the collection of edge-extending lines
 $\{ l[e],\; e \in \Edges(\gamma) \}$ coincides with the collection of lines
 determined by $\Lambda^{\cal M}_D.$
 This shows that the field traced by the particles coincides in law 
 with ${\cal A}^{\cal M}_D$ and that (\ref{ZZ}) holds, thus completing
 the proof of the theorem. $\Box$ 

\subsection{Conclusions from the dynamic representation}\label{CDYNRES}
 The crucial property of polygonal fields ${\cal A}^{\cal M}_{(\cdot)}$ for
 $\beta = 1$ is their {\it consistency}  stated in Theorem \ref{WlasnosciZgd} below
 and established using the dynamic representation from Theorem \ref{ReprDyn} in analogy
 with  Arak \& Surgailis \cite{AS1}. This is where the name {\it consistent regime} comes
 from. Another important feature of ${\cal A}^{\cal M}_{(\cdot)}$ is the explicit
 knowledge of its linear sections, see Theorem \ref{WlasnosciZgd}, which can be
 further used to obtain information about the first and second order correlation
 structure of the field, see Theorem \ref{Corr12} and Corollaries
 \ref{FOD}, \ref{SOD} and \ref{PPcorrC}.

 \begin{theorem}\label{WlasnosciZgd}
  The polygonal field ${\cal A}^{\cal M}_D$ enjoys the following properties
 \begin{description}
  \item{\bf Consistency:} For bounded open convex $D' \subseteq D \subseteq {\Bbb R}^2$
        the field ${\cal A}^{\cal M}_D \cap D'$ coincides in law with
        ${\cal A}^{\cal M}_{D'}.$ This
        allows us to construct the whole
        plane extension of the process ${\cal A}^{\cal M}$ such that
        for each bounded open convex $D \subseteq {\Bbb R}^2$ the field
        ${\cal A}^{\cal M}_D$ coincides in law with 
        ${\cal A}^{\cal M} \cap D.$ 
  \item{\bf Linear sections:} For a straight line $l$ in ${\Bbb R}^2$ the collection of
        intersection points and intersection directions of $l$ with the edges of the
        polygonal field ${\cal A}^{\cal M}$ coincide in law with the
        corresponding collection for the Poisson line process $\Lambda^{\cal M}.$ 
 \end{description}
 \end{theorem}

 \paragraph{Proof}
  To establish the {\bf Consistency} property, choose a bounded open convex set
  $D \subseteq {\Bbb R}^2$ and a straight line $l$ intersecting $D$ and define
  $D'$ to be the set of points of $D$ lying to the left from $l$ (lower time coordinates).
  Clearly then, from the dynamic representation Theorem \ref{ReprDyn} we conclude the
  {\bf Consistency} statement for the so chosen $D$ and $D'.$ Noting that the dynamic
  representation is equally available upon rotating the space-time coordinate system
  we see that the {\bf Consistency} holds as well upon cutting off the part of the set
  $D$ lying to the left  of $l.$ This means however that the consistency holds upon
  cutting off pieces of the original set with arbitrary straight lines \--- a repetitive
  use of this procedure and a possible passage to the limit allows us to {\it carve} from $D$
  its arbitrary convex subset. This proves the {\bf Consistency} claim. To establish the
  {\bf Linear sections} statement pick a straight line $l$ and choose the space-time
  coordinate system so that $l$ coincides with its spatial axis. The {\bf Linear sections}
  claim follows now from the form of the boundary birth mechanism in the dynamic
  representation in view of the {\bf Consistency} property. $\Box$

\section{Generalised dynamic representation for consistent fields}\label{GDR}
 The dynamic construction of consistent polygonal fields, originating in the homogeneous set-up
 from Arak \& Surgailis \cite{AS1} and discussed in Subsection \ref{CONSIS} above, can be regarded
 as {\it revealing} increasing portions of the polygonal field in the course of the time flow.
 Under this interpretation, the portion of a polygonal field in a bounded open convex domain
 $D$ {\it uncovered} by time $t$ is precisely its intersection with 
 $D_t = \bar{D} \cap (-\infty,t] \times {\Bbb R}_+.$  
 The idea underlying our generalised version of dynamic representation developed in the present
 section is to replace the above family $D_t$ by some other time-increasing family of
 subsets of $D$, also denoted $D_t$ in the sequel, eventually covering the whole $D,$ and to
 try to provide a natural construction of the polygonal field being gradually uncovered on
 the growing domain $D_t$ in the course of the time flow. We shall always assume that
 $D_t$ be convex for otherwise we would have to deal with situations where two or more
 disconnected parts of an edge of the field have been revealed which leads to unwanted
 and cumbersome dependencies along the segments connecting these parts. Taking this
 into account and having formal convenience in mind we impose the following natural
 assumptions on $D$ and $D_t, \; t\in [0,1],$
 \begin{description}
  \item{\bf (D1)} $(D_t)_{t \in [0,1]}$ is a strictly increasing family of compact
        convex subsets of $\bar D = D \cup \partial D.$
  \item{\bf (D2)} $D_0$ is a single point $x$ in $\bar{D} = D \cup \partial D.$
  \item{\bf (D3)} $D_1$ coincides with $\bar D.$
  \item{\bf (D4)} The family $(D_t)_{t \in [0,1]}$ enjoys the property that 
         ${\cal M}(\{ l,\; \exists_{t \in [0,1]} \card (l \cap \partial D_t) > 2 \}) = 0,$
         in particular $\bar{D}$ itself satisfies 
         ${\cal M}(\{ l,\; \card (l \cap \partial D) > 2 \}) = 0.$
  \item{\bf (D5)} $D_t$ is continuous in the usual Hausdorff metric on compacts.  
 \end{description}
 Under these conditions, for ${\cal M}$-almost each $l \in [[D]]$ the intersection
 $l \cap D_{\tau_l}$ consists of precisely one point ${\Bbb A}(l),$ where
 $\tau_l = \inf \{ t \in [0,1],\; D_t \cap l \neq \emptyset \}.$ 
 The point ${\Bbb A}(l)$ will be referred to as the {\it anchor point} for $l,$
 this induces the {\it anchor} mapping ${\Bbb A}: [[D]] \to D.$
 Consider now the following dynamics in time $t \in [0,1],$ with all updates
 given by the rules below performed independently of each other.
 \begin{description}
  \item{\bf (GE:Initialise)} Begin with empty field at time $0,$
  \item{\bf (GE1)} Between critical moments listed below, during the time interval $[t,t+dt]$ the
                   unfolding field edges in $D_t$ reaching $\partial D_t$ extend straight to
                   $D_{t+dt} \setminus D_t,$
  \item{\bf (GE2)} When a field edge hits the boundary $\partial D,$ it stops growing in this direction
                  (recall that ${\cal M}$-almost everywhere the intersection of a line with 
                    $\partial D$ consists of at most two points),
  \item{\bf (GE3)} When two unfolding field edges intersect in $D_{t+dt} \setminus D_t,$ they are not 
                   extended any further beyond the intersection point (stop growing in the 
                   direction marked by the interesction point),
  \item{\bf (GE4)} A field edge extending along $l \in [[D_t]]$ updates its direction during
                   $[t,t+dt]$ and starts unfolding along $l' \in [[l^{[t,t+dt]}]],$
                   extending away from the anchor point ${\Bbb A}(l'),$ with probability
                   ${\cal M}(dl'),$ where $l^{[t,t+dt]} := l \cap (D_{t+dt} \setminus D_t).$
                   Directional updates of this type are all performed independently,
  \item{\bf (GE:LineBirth)} Whenever the anchor point ${\Bbb A}(l)$ of a line $l$ falls
                  into $D_{t+dt} \setminus D_t,$ the line $l$ is born at the time $t$ at
                  its anchor point with probability ${\cal M}(dl),$ whereupon it begins
                  extending in both directions with the growth of $D_t$ (recall that
                  $l$ is ${\cal M}$-almost always tangential to $\partial D_t$ here),
  \item{\bf (GE:VertexBirth)} For each intersection point of lines $l_1$ and $l_2$ falling
                  into $D_{t+dt} \setminus D_t,$ the pair of field lines $l_1$ and $l_2$ is
                  born at $l_1 \cap l_2$ with probability ${\cal M}(dl_1) {\cal M}(dl_2),$
                  whereupon both lines begin unfolding in the directions away from their
                  respective anchor points ${\Bbb A}(l_1)$ and ${\Bbb A}(l_2).$ 
 \end{description}
 Observe that the evolution rule {\bf (GE:VertexBirth)} means simply that pairs of
 lines are born at birth sites distributed according to a Poisson point process in
 $D$ with intensity measure $\lnaw {\cal M} \pnaw,$ in analogy to the standard
 dynamic representation in Subsection \ref{DYNRES}. The essential difference
 though is that here the pairs of lines emitted from a birth site extend away from
 their respective anchor points rather than always in the direction determined by a
 single time axis. Further, the line birth events given by {\bf (GE:LineBirth)} replace
 the boundary birth events in the standard dynamic representation.   
 It is also worth noting that if we choose the family $D_t$ so that
 $D_t := \bar D \cap (-\infty,(1-t)x_{\min} + t x_{\max}] \times {\Bbb R}_+,$
 where $x_{\min}$ and $x_{\max}$ are the minimal and maximal x-coordinates of a
 point in $D$ (assume that $\bar D$ contains exactly one point with x-coordinate
 $x_{\min}$ and that ${\cal M}$ assigns zero mass to the set of vertical lines),
 the generalised dynamic representation {\bf (GE)} coincides with
 the standard Arak \& Surgailis one determined by rules {\bf (E1-4)} in Subsection
 \ref{DYNRES} above, and we have ${\Bbb A}(l) = \iin(l,D).$

 In analogy with the corresponding result for the usual dynamic construction, as
 established in the proof of Theorem \ref{ReprDyn}, we show that the field
 resulting from the above {\bf (GE)} construction does coincide in law with
 ${\cal A}^{\cal M}_D.$ 
 \begin{theorem}\label{GEthm}
  The random contour ensemble resulting from the above construction {\bf (GE)}
  coincides in law with ${\cal A}^{\cal M}_D.$
 \end{theorem}

 \paragraph{Proof}
  The proof is very similar to that of Theorem \ref{ReprDyn}. We pick some $\gamma \in \Gamma_D$
  and calculate the probability that the outcome of the above dynamic construction falls
  into $d\gamma.$ To this end, we note that:
  \begin{itemize}
   \item Each edge $e \in \Edges(\gamma)$ containing the anchor point ${\Bbb A}(l[e])$
         and hence resulting from a line birth event due to the rule {\bf (GE:LineBirth)},
         contributes to the considered probability the factor  ${\cal M}(dl[e])$ (line
         birth probability for $l[e]$) times $\exp(-{\cal M}([[e]]))$ (no directional
         updates along $e$),
   \item Each of the two edges $e_1,e_2 \in \Edges(\gamma)$ steming from a common interior
         birth vertex $l[e_1] \cap l[e_2]$ yields the factor ${\cal M}(dl[e_i]),\;
         i=1,2,$ (coming from the vertex birth probability due to the rule 
         {\bf (GE:VertexBirth)}) times $\exp(-{\cal M}([[e_i]]))$ (no directional updates along $e_i$),
  \item  Each of the edges $e \in \Edges(\gamma)$ arising due to a directional update 
         in {\bf (GE4)} yields the factor ${\cal M}(dl[e])$ (directional update probability)
         times $\exp(-{\cal M}(dl[e]))$ (no directional updates along $e$),
  \item  The absence of interior birth sites in $D \setminus \gamma$ yields the factor of
         $\exp(-\lnaw {\cal M} \pnaw(D)),$
  \item  Finally, the absence of line birth events for all lines in $[[D]]$ except for
         the finite collection $\{ l[e],\; e \in \Edges(\gamma) \}$ yields
         the additional factor $\exp(-{\cal M}(D)).$
 \end{itemize}      
 Putting these observations together we conclude that the probability element of
 $\gamma$ resulting from the generalised construction above is 
 $$ \frac{\exp\left(-L^{\cal M}(\gamma)\right)}
         {\exp(\lnaw{\cal M}\pnaw(D))}
    \frac{\prod_{e\in\Edges(\gamma)} {\cal M}(dl[e])}
         {\exp({\cal M}(D))} = 
    \frac{\exp\left(-L^{\cal M}(\gamma)\right)}
         {\exp(\lnaw{\cal M}\pnaw(D))}        
         d{\Bbb P}(\gamma \in \Gamma_D(\Lambda^{\cal M})) $$
 and thus the field obtained by this construction coincides in law 
 with ${\cal A}^{\cal M}_D$ as required. This completes
 the proof of the theorem. $\Box$

 \section{Correlations in the consistent regime}\label{CORCR}
 In the present section we use the dynamic construction and its generalised
 version to describe the correlation structure of the consistent field
 ${\cal A}^{\cal M}.$ Due to the polygonal nature of the considered field
 the natural object of our interest are the {\it edge correlations}
 \begin{equation}\label{CORR}
  \sigma^{\cal M}[dl_1,x_1;\ldots;dl_k,x_k] := {\Bbb P}\left( \forall_{i=1}^k 
  \exists_{e \in \Edges({\cal A}^{\cal M})} \; \pi_{l_i}(x_i) \in e,\; l[e] \in dl_i \right),
 \end{equation}
 where $l_1,\ldots,l_k$ are straight lines and $\pi_{l_i}$ is the orthogonal
 projection on $l_i.$ In almost all cases below we shall be interested in correlations
 with $x_i \in l_i,$ in which case $\sigma^{\cal M}[dl_1,x_1;\ldots;$ $dl_k,x_k]$ 
 can be interpreted as the probability element that the polygonal field ${\cal A}^{\cal M}$
 passes through points $x_i$ in the directions determined by the respective lines
 $l_i,\; i=1,\ldots,k.$ For general $x_i,$ not necessarily lying on $l_i,$ the $k$-fold
 correlation $\sigma^{\cal M}[dl_1,x_1;\ldots;dl_k,x_k]$ is the probability that
 the polygonal field passes through points $\pi_{l_i}(x_i)$ in the directions
 determined by the respective lines $l_i,\; i=1,\ldots,k.$  
 
\subsection{First and second order edge correlations} 
 The first and second order edge correlations are easily determined using 
 the {\bf Linear sections} statement of Theorem \ref{WlasnosciZgd}
 concluded from the standard dynamic representation. 
 \begin{theorem}\label{Corr12}
  We have for $x \in l$
  $$ \sigma^{\cal M}[dl,x] = {\cal M}(dl) $$
  and, for $x_1 \in l_1,\; x_2 \in l_2,$
  $$ \sigma^{\cal M}[dl_1,x_1;dl_2,x_2] =
     \left\{ \begin{array}{ll} {\cal M}(dl_1) {\cal M}(dl_2), & \mbox{ if } l_1 \neq l_2, \\
             \exp(-2 {\cal M}([[ \overline{x_1 x_2}]])) {\cal M}(dl), & \mbox{ if } l_1 = l_2, 
     \end{array} \right. $$
  where $\overline{x_1 x_2}$ is the segment joining $x_1$ to $x_2.$ 
 \end{theorem}

 \paragraph{Proof}
  The statement for $\sigma^{\cal M}[dl,x]$ and $\sigma^{\cal M}[dl_1,x_1;dl_2,x_2],\;
  l_1 \neq l_2,$ is
  a direct consequence of the {\bf Linear sections} part of Theorem \ref{WlasnosciZgd}
  and of the properties of the Poisson line process $\Lambda^{\cal M}.$ To find 
  $\sigma^{\cal M}[dl,x_1;dl,x_2]$ for $x_1, x_2 \in l$ note that in order to 
  have an edge $e$ of the field passing through both $x_1$ and $x_2$ along $l,$  
  we have to ensure that
   \begin{itemize}
    \item There is $e \in\Edges({\cal A}^{\cal M})$ such that $l[e] \in dl$
          and $x_1,$ which happens with probability ${\cal M}(dl)$ in view of
          {\bf Linear sections} in Theorem \ref{WlasnosciZgd},
    \item There are no other edges of the field crossing $e$ between $x_1$ and $x_2,$
          which happens with probability $\exp(-2 {\cal M}([[\overline{x_1 x_2}]])).$
          Indeed, to see it choose the spatial axis in the standard dynamic representation in
          Subsection \ref{DYNRES} very close to $l$ (not exactly $l$ to keep the
          velocity of the particle tracing $e$ finite though very large)
          and observe that $e$ can be crossed by
          \begin{itemize}
           \item space-time trajectories of particles coming from the past, which 
                 happens with probability $1 - \exp(-{\cal M}([[\overline{x_1 x_2}]]))[1+o(1)]$
                 by the {\bf Linear sections} property,
           \item edges arising due to velocity updates along $\overline{x_1 x_2},$
                 which happens with 
                 probability $1 - \exp(-{\cal M}([[\overline{x_1 x_2}]]))[1+o(1)]$
                 by the dynamic rule {\bf (E4)}.  
          \end{itemize}
          Since the velocity update events in {\bf (E4)} are independent of the past
          particle configuration, letting the spatial axis approach $l$ we obtain the
          required conclusion.
    \end{itemize}
  Combining the factors listed above we obtain the required formula for
  $\sigma^{\cal M}[dl,x_1;dl,x_2],$ thus completing the proof of the Theorem. $\Box$\\

 \subsection{Higher order edge correlations for general activity measures}\label{HOCgen}
  To describe the higher order correlation structure of the field ${\cal A}^{\cal M}$
  we need the full power of the more flexible generalised dynamic representation. 
  Consider a collection $(l_1,x_1),(l_2,x_2),$ $\ldots,(l_k,x_k),\; k \geq 1,$ of 
  pairwise different lines $l_i$ and points $x_i$ with $x_i \in l_i$ and
  $x_i \not\in l_j$ for $j \neq i.$ Such collections are said to be in
  {\it general position} below and such assumption will be imposed
  an all collections $(l_i,x_i)$ considered in this subsection, often
  without a further mention. Also, throughout this subsection we always
  assume for formal convenience that all $x_i,\; i=1,\ldots,n,$ are contained
  in a bounded open convex set $D,$ playing the usual role of the field
  domain. Clearly, the particular choice of $D$ is irrelevant due to the
  consistency of the field ${\cal A}^{\cal M}.$
  
  We say that the edge correlations of the field ${\cal A}^{\cal M}$ 
  {\it factorise} on a collection  $(l_i,x_i)_{i=1}^k$ if
  $\sigma^{\cal M}[dl_1,x_1;\ldots;dl_k,x_k]$ coincides with the
  product $\prod_{i=1}^k \sigma^{\cal M}[dl_i,x_i].$ Moreover, we
  say that the {\it collective factorisation of correlations} holds
  for a family  $(l_i,x_i)_{i=1}^k$ if the correlations factorise for
  $(l_i,x_i)_{i=1}^k$ and all its sub-collections. In view of Theorem \ref{Corr12}
  $(l_i,x_i)_{i=1}^k$ collectively factorises iff
  \begin{equation}\label{CollFct}
     \sigma^{\cal M}[dl_{i_1},x_{i_1};\ldots;dl_{i_m},x_{i_m}]
     = \prod_{j=1}^m {\cal M}(dl_{i_j})
  \end{equation}
  for all sub-collections $(l_{i_j},x_{i_j})_{i=1}^m.$
 
 \paragraph{General collective factorisation problem and precedence graphs}
  Our main objective in the
  present subsection is to find general conditions characterising
  collections $(\bar l,\bar x)$ in general position enjoying the collective factorisation
  property. To see what may be plausible answers to this question let us make
  the following basic observations very helpful in interpreting the geometry
  of the higher order edge correlations. For a collection
  $(\bar l,\bar x) = (l_1,x_1;\ldots;l_k,x_k),\; x_i \in l_i,$ in general position by 
  $\Gamma(\bar l,\bar x)$ we shall mean the family of admissible polygonal configurations
  $\gamma$ in the plane consisting of $k$ edges $e_1,\ldots,e_k$ such that 
  $e_i$ lies on $l_i$ and contains $x_i$ for all $i=1,\ldots,k.$ Looking at
  the structure of the family $\Gamma(\bar l,\bar x)$ is closely related to the 
  following natural problem which marks its presence in various
  domains ranging from kid games to studies on human and computer vision,
  see \cite{AC07} and the references therein: given the collection
  $(\bar l,\bar x) = (l_i,x_i)_{i=1}^k$
  draw a family of closed curves (here a polygonal configuration) such that 
  each point $x_i$ lies on one of the curves and the direction of the
  curve at $x_i$ is determined by $l_i.$ Clearly, the solution to this
  problem is non-unique in general, yet if we require in addition that
  the configuration we draw belong to $\Gamma(\bar l,\bar x),$ that is to say we
  may only draw over the lines $l_i$ and on each of these lines we
  have to draw precisely one segment of non-zero length, then it is
  often the case that $\Gamma(\bar l,\bar x)$ is a singleton. Now, assume that
  all points $x_i$ as well as the intersection points $y_{i,j}$
  for $l_i$ and $l_j$ lie very close to each other, say they are all
  contained in a disk of very small radius $r,$ then it follows from
  (\ref{GIBBS1}) and (\ref{GIBBS2}) combined with the consistency property
  of the field that 
  \begin{equation}\label{MALER}
   \sigma^{\cal M}[dl_1,x_1;\ldots;dl_k,x_k] =
    N(\bar l,\bar x) {\cal M}(dl_1) \ldots {\cal M}(dl_k) (1+o_r(1)),
   \end{equation}
   where $N(\bar l,\bar x)$ is the cardinality of $\Gamma(\bar l,\bar x).$
   Indeed, to see it note that the Boltzmann weight $\exp(-L^{\cal M}(\cdot))$ tends to
   $1$ as $r \to 0.$  Thus, the 
   factorisation holds in small $r$ asymptotics precisely when
   $\Gamma(\bar l,\bar x)$ is a singleton. As will be discussed in the
   sequel, there are many collections for which $N(\bar l,\bar x) = 0$ or
   $N(\bar l,\bar x) > 1.$ 
   In particular, it is not the case that the factorisation holds for
   all collections in general position.
   Neither can it be hoped though 
   that the formula (\ref{MALER}) holds in general non-asymptotic regime:
   as we shall show in Theorem \ref{EDecT} in Section \ref{CORRDEC} in the
   particular case of rectangular fields, when the distances between $x_i$'s
   get large, the correlations converge exponentially fast to the product,
   that is to say the polygonal field exhibits exponential decay of
   dependencies. To conclude these considerations, we say that $(l_i,x_i)_{i=1}^k$ 
   enjoys collective factorisation on all scales iff 
   $(\alpha l_i,\alpha x_i)_{i=1}^k$ factorises collectively for each
   $\alpha > 0,$ where by $(\alpha l, \alpha x)$ we understand the
   re-scaled version of $(l,x)$ with scaling factor $\alpha.$ 
   Then the above discussion shows that    
   \begin{lemma}\label{Rysowanie}
    For a collection $(\bar l,\bar x) = (l_i,x_i)_{i=1}^k$ in general position
    a necessary condition to enjoy
    collective factorisation on all scales is that $N(\bar l,\bar x) = 1$ and
    $N(\bar l',\bar x') = 1$ for all non-empty sub-collections $(\bar l',\bar x')$
    of $(\bar l,\bar x).$
   \end{lemma}
   To reformulate this condition in more tangible terms, we build for      
   each family $(l_i,x_i)_{i=1}^k$ its {\it precedence graph} 
   ${\cal G}[l_1,x_1;\ldots;l_k,x_k]$ as follows. We split each $l_i$
   at $x_i$ into two half-lines both oriented in the directions away from
   $x_i.$ This creates a directed graph ${\cal G}[l_1,x_1;\ldots;l_k,x_k]$
   whose vertices are the {\it generating points} $x_i$ and the
   {\it intersection points} $y_{i,j}$ of respective pairs of lines
   $l_i, l_j.$ As we shall see below, the name {\it precedence graph} comes
   from its relationship to the order in which the points of $\bigcup_{i=1}^k l_i$
   are revealed in the course of a suitable instance of the generalised
   graphical construction. For now, we claim that
   \begin{lemma}\label{Bezcyklowosc}
    For a collection $(\bar l,\bar x) = (l_i,x_i)_{i=1}^k$ in general position
    the following conditions are equivalent
    \begin{enumerate}
     \item $N(\bar l,\bar x) = 1$ and $N(\bar l',\bar x') = 1$ for all non-empty
           sub-collections of $(\bar l,\bar x),$
     \item The precedence graph ${\cal G}[l_1,x_1;\ldots;l_k,x_k]$ is acyclic.
    \end{enumerate}
   \end{lemma}

   \paragraph{Proof} To show that the first condition implies the second one
    note that if the precedence graph contains a cycle of length $m$ built by
    $(l_{i_1},x_{i_1};\ldots;l_{i_m},x_{i_m}),$ then for the sub-collection
    $(\bar l',\bar x') = (l_{i_j},x_{i_j})_{j=1}^m$ we have $N(\bar l',\bar x') = 0$ if
    $m$ is odd and $N(\bar l',\bar x') = 2$ if $m$ is even.   
   
    To prove the inverse implication assume that the precedence graph 
    ${\cal G}[l_1,x_1;\ldots;l_k,x_k]$ contains
    no cycles and observe that this induces a partial ordering on
    $\bigcup_{i=1}^k l_i$ in which the $x_i$'s are minimal points, each $x_i$
    first on its corresponding line $l_i$ and with the remaining points on
    $l_i$ ordered according to the natural linear orderings directed away from
    $x_i$ on the two half-lines. This will be referred to as the {\it structural
    ordering} for the collection $(l_i,x_i)_{i=1}^k$ in the sequel. Consider the
    following incremental construction of a graph belonging to $\Gamma(\bar l,\bar x),$
    where at each step we obtain a collection of $k$ segments
    $\iota_i \subseteq l_i,\; i=1,\ldots,k$ with non-intersecting interiors
    but possibly sharing endpoints and possibly degenerated to $x_i$'s,
    eventually yielding the entire graph under construction.  
    \begin{enumerate}
     \item begin with  $\iota_i := \{ x_i \},\; i=1,\ldots,k,$
     \item choose an intersection point $y_{i,j}$ which: 
      \begin{enumerate}
       \item belongs to neither $\iota_i$ nor $\iota_j$
             and hence to neither of the remaining
             segments $\iota_l.$
       \item enjoys the property that extending both $\iota_i$
             and $\iota_j$ to contain $y_{i,j}$ creates no
             T-shaped or X-shaped nodes (vertices of order
             three or four respectively), which is equivalent
             to both $\iota_i$ and $\iota_j$ having their endpoints
             pointing towards $y_{i,j}$ not shared with any other segment
             $\iota_l$. (Such endpoints will be referred to as {\it loose
             ends} below. Note that by definition of our construction
             procedure the only possible loose end on $\iota_i$ is $x_i.$)  
       \item is minimal with the above two properties.
      \end{enumerate}
       with possible ties broken in an arbitrary way,  
     \item extend both segments $\iota_i$ and $\iota_j$ to contain
           $y_{i,j}$ \--- this action will be referred to as
           {\it adding} $y_{i,j}$ {\it to the graph} for short in
           the sequel of this argument,  
     \item return to 2. unless no more $y_{i,j}$'s can be added,
     \item whenever a segment $\iota_l$ of the constructed graph ends with
           a loose node, that is to say a vertex of order one, extend this
           segment to the half-line in the direction of the node.
    \end{enumerate}
    In other words, we initialise our graph under construction with the set of
    generating points $\{x_1,\ldots,x_k\}$ whereupon we let it grow along the lines
    $l_i,$ adding subsequent intersection points $y_{i,j}$ in the order
    determined by the structural ordering whenever this does not violate
    the usual constraints imposed on a polygonal configuration, and discarding
    those $y_{i,j}$'s whose addition would violate these constraints. 
    It is easily seen that this procedure yields in a finite number of steps
    a graph belonging to $\Gamma(\bar l,\bar x),$ thus in particular $\Gamma(\bar l,\bar x)
    \neq \emptyset$ and $N(\bar l,\bar x) \geq 1.$ To show that $N(\bar l,\bar x) \leq 1$
    note that  in fact {\it all} graphs in $\Gamma(\bar l,\bar x)$ can be obtained by the
    above procedure. Indeed, this is seen inductively. First, all $x_i$'s  have
    to belong to such a graph. Second, whenever a graph in $\Gamma(\bar l,\bar x)$ contains
    some configuration of $\iota_i$'s as its (sub)segments and some of these $\iota_i$'s
    have loose ends $x_i$'s, each of such $\iota_i$'s has to be extended, either up to
    intersection with the extension of another segment or to the entire half-line if such
    an intersection is not available. This is because no loose ends can be present in the final
    graph and because we cannot hope that the potential loose end $x_i$ on $\iota_i$
    could be possibly reached by another edge $\iota_j,\; j \neq i,$ since $x_i \not\in l_j,\;
    j \neq i.$ This means that at least one intersection point $y_{i,j}$ as in 2.(a,b,c)
    has to be added to the graph as in 2. Consequently, two different
    graphs in $\Gamma(\bar l,\bar x)$ could only arise in the course
    of two distinct instances of the above procedure, with different choices of minimal
    points $y_{i,j}$ in step 2. We argue that the so constructed graphs necessarily
    coincide though. To see it let us make the following observations
    \begin{itemize}
     \item If at a certain stage of the construction some intersection point 
           $y_{i,j}$ becomes admissible for 2.(a,b,c), there is no way to make
           it inadmissible whatever choices be made in the sequel of the
           construction save for adding $y_{i,j}$ to the graph. Indeed, 
           to make $y_{i,j}$ inadmissible in the sequel of the construction
           we would have to close either $\iota_i$ or $\iota_j$ by intersecting
           its extension with an extension of some other $\iota_l$ before
           $\iota_i$ or $\iota_j$ hits $y_{i,j},$ which is not possible by
           the minimality of $y_{i,j}$ in 2.(a,b,c).   
     \item If at a certain stage of the construction some intersection point 
           $y_{i,j}$ becomes admissible for 2.(a,b,c), there is no way to
           make it inadmissible by adding some other intersection points in
           the prequel of the construction in agreement with the rule 2.
           (without removing those already present though). Indeed, to make
           $y_{i,j}$ inadmissible in this way we would have to have added some
           other intersection point $y$ on $\iota_i$ or $\iota_j$ prior to
           the considered construction stage, but if this were possible this
           point $y$ would have to reach its admissibility before $y_{i,j}$ and
           in view the previous observation there would be no way to make
           it inadmissible before $y_{i,j}$ reaches its admissibility,
           which would contradict the minimality of $y_{i,j}$ in 2.(a,b,c)
           at the moment of becoming admissible.
    \end{itemize}
   Consequently, regardless of any particular sequence of choices in a given
   instance of our construction, the points admissible for 2.(a,b,c) in the
   very first step of the construction (call them the {\it first generation})
   eventually have to be added in view of the first observation above. Next,
   by the first and second observation above, all points which become
   admissible upon adding just the first generation
   (call them the {\it second generation}) also have to be added at some stage
   of the construction regardless of the particular sequence of choices made.
   Proceeding further this way we conclude inductively that the so-defined
   generations of all orders will eventually be added to the graph. 
   However such generation-wise addition of intersection points (plus
   possible extensions to half-lines at the end of the procedure) is also
   an instance of our incremental construction and it yields a valid graph
   belonging to $\Gamma(\bar l,\bar x)$ which, in view of the above discussion,
   enjoys the property of being contained in any other graph obtained in
   any instance of our construction. It remains to observe that a valid 
   graph in $\Gamma(\bar l,\bar x)$ cannot be further extended by our 
   construction. All this means that the results of all possible instances
   of our construction coincide and the order in which the points $y_{i,j}$
   are considered in 2. is irrelevant. This way, we have shown that
    $N(\bar l,\bar x) = 1$ as required. $\Box$

   We ask whether the necessary precedence graph acyclicity condition stated
   in combined Lemmas \ref{Rysowanie} and \ref{Bezcyklowosc} is also a sufficient
   condition for collective 
   factorisation of correlations (on all scales). Below we are going to show
   that the answer to this question is positive for rectangular fields, see
   Subsection \ref{FUCH} and Theorem \ref{RectangFactor} there.
   For general fields we were only able to establish
   collective edge factorisation under a somewhat stronger sufficient condition
   though, see Theorem \ref{CORRH}, and we do not know at present if this condition
   can be weakened. 
  \paragraph{General sufficient condition for collective factorisation}
   We proceed with the general case first. To this end, we note that the
   generalised dynamic construction with its {\bf (GE:LineBirth)} rule
   allows us to conclude:
   \begin{lemma}\label{CORRHC}
     Assume that $(l_i,x_i)_{i=1}^k$ is a collection of pairwise distinct
     lines and points in general position and such that $l_1$ does not hit
     the convex hull $\conv(\{x_2,\ldots,x_k\}).$ Then
     $$ \sigma^{\cal M}[dl_1,x_1;\ldots;dl_k,x_k] = {\cal M}(dl_1)
        \sigma^{\cal M}[dl_2,x_2;\ldots;dl_k,x_k]. $$
    \end{lemma} 
  Indeed, under the assumptions of the corollary, by standard geometry it is
  always possible to construct an increasing family $(D_t),\; t \in [0,1],$ of
  convex  compacts satisfying rules {\bf (D1-5)} as in Section \ref{GDR}
  and such that $D_t$ covers the whole of $\conv(\{x_2,\ldots,x_k\})$ before
  hitting $l_1$ and that ${\Bbb A}(l_1) = x_1,$ that is to say $x_1$ is the
  first point of $l_1$
  hit by $D_t.$ This yields the required statement by the consistency
  of the field and by the {\bf (GE:LineBirth)} rule. This way of thinking
  suggests a natural sufficient condition for collective factorisation to
  hold under general activity measures,
  namely that there exist an increasing family $(D_t),\; t \in [0,1],$ of convex
  compacts as in Section \ref{GDR} and such that $x_i = {\Bbb A}(l_i)$ for all
  $1 \leq i \leq k,$ that is to say $x_i$ is the point at which $l_i$ is first
  hit by $D_t.$ However, the so formulated collective factorisation condition
  is rather untractable, therefore we look for its equivalent reformulation
  in more tangible terms. To this end we {\it augment} the precedence graph
  ${\cal G}[l_1,x_1;\ldots;l_k,x_k]$ as follows. For each two lines $l_i, l_j$
  intersecting at some $y_{i,j}$ we note that $l_i$ and $l_j$ divide the plane
  into four regions, one with both $x_i$ and $x_j$ on its boundary,
  two with either $x_i$ or $x_j,$
  and finally one, termed the {\it trap region}, with neither. 
  Now, if some $x_m$ falls into such a trap region, it is easily seen that for
  each family $D_t$ such that ${\Bbb A}(l_i) = x_i,\; i=1,\ldots,k,$ we must
  have $y_{ij}$ hit by $D_t$ prior to $x_m.$ Indeed, if at a certain time $D_t$
  contains $x_m,$ then if by that time it contains either of the points $x_i$
  or $x_j,$ it has to intersect both $l_i$ and $l_j$ and thus contain both
  $x_i, x_j$  and hence also $y_{ij}$ by convexity. On the other
  hand if $D_t$ hits $x_m$ before hitting any of $x_i$ and $x_j$ then we cannot
  have simultaneously ${\Bbb A}(l_i) = x_i$ and  ${\Bbb A}(l_j) = x_j.$
  Taking this into account we add in ${\cal G}[l_1,x_1;\ldots;l_k,x_k]$
  a directed {\it trap edge} from $y_{ij}$ to $x_m$ for each $i,j,m$ as above.
  Denote the resulting {\it augmented} directed graph by ${\cal G}^+[l_1,x_1;\ldots;l_k,x_k].$
  It is important to observe that, as follows by its construction and the above
  discussion, the orientation of edges in this graph indicates
  the order in which its vertices are hit by the sought for increasing family
  $D_t,$ should it exist. Consequently, the acyclicity of the augmented precedence
  graph is a necessary condition for the existence of such $D_t.$ In the proof 
  of the following crucial theorem we show that it is also a sufficient condition. 
  \begin{theorem}\label{CORRH}
   Assume that $(l_i,x_i)_{i=1}^k$ is a collection of pairwise distinct
   lines and points in general position such that the augmented precedence graph 
   ${\cal G}^+[l_1,x_1;\ldots;l_k,x_k]$ is acyclic. Then $(l_i,x_i)_{i=1}^k$ admits
   collective factorisation of correlations on all scales. 
  \end{theorem}

  \paragraph{Proof} In context of the discussion above and in view of the generalised
    graphical construction and its {\bf (GE:LineBirth)} rule as used in 
    Lemma \ref{CORRHC}, it suffices to show that there exists an increasing
    family $D_t$ such that ${\Bbb A}(l_i) = x_i$ for all $i=1,\ldots,k.$ Now,
    to establish the existence of such a family it is enough to know that
    having  ${\cal G}^+[l_1,x_1;\ldots;l_k,x_k]$ acyclic implies the existence of a
    permutation $(s_i)_{i=1}^k$ of indices such that $l_{s_{i+1}}$ does
    not hit the convex hull $\conv(\{x_{s_1},\ldots,x_{s_i}\}),$ $i=1,\ldots,k-1.$
    Indeed, the family $D_t$ is then easily constructed by induction in $k:$ 
    it starts growing from $x_{s_1}$ whereupon it hits the consecutive
    points $x_{s_2},x_{s_3},\ldots,$ and since $l_{s_{k+1}}$ is disjoint
    with $\conv(\{ x_{s_1},\ldots,x_{s_k}\}),$ the family $D_t$ can be
    chosen so that it does not hit $l_{s_{k+1}}$ until it reaches all previous
    lines $l_{s_i},\; i \leq k,$ and so that ${\Bbb A}(l_{s_i}) = x_{s_i}$ for
    $i \leq k$ (by induction hypothesis) and then ${\Bbb A}(l_{s_{k+1}}) = x_{s_{k+1}}$
    (again by disjointness of $l_{s_{k+1}}$ with the convex hull of preceding
     $x_{s_i}$'s). This can be equivalently interpreted as iterative application
    of Lemma \ref{CORRHC}.    

    To proceed, we use induction in $k$ to show that whenever ${\cal G}^+[\ldots]$
    is acyclic, the required permutation ensuring the disjointness of lines with
    convex hulls of sets of preceding points exists. To this end, observe first
    that our statement trivialises for $k=1.$ Next, choose a collection
    $(l_1,x_1;\ldots;l_k,x_k)$ and note that if there exists a line in the
    collection which does not hit the convex hull generated by all the remaining
    points then by Lemma \ref{CORRHC} we can resort to the inductive hypothesis for the
    collection with the line removed. Consequently, we only have to show that if the
    collection $(l_1,x_1;\ldots;l_k,x_k)$ is such that {\it each} line $l_i$ hits the
    convex  hull $\conv(\{ x_j,\; j \neq i \})$ then ${\cal G}^+[l_1,x_1;\ldots;l_k,x_k]$
    contains a cycle. Again resorting to inductive argument if needed, we may
    assume without loss of generality that $(l_1,x_1;\ldots;l_k,x_k)$ is minimal
    with this property, that is to say it does not contain a proper subcollection
    such that each its line hits the convex hull generated by the remaining points.
    This minimality assumption implies that all $x_i$'s are extreme points 
    (vertices) of the convex hull $\conv(\{x_1,\ldots,x_k\}),$ for otherwise
    removing a non-extreme point we would obtain a subcollection enjoying the
    considered property. For formal convenience we let the vertices $x_i$
    be ordered clockwise along the boundary of the convex hull and we
    interpret the indices modulo $k,$ that is to say $x_{k+1} = x_1$ etc.
    By minimality, should we remove a vertex $x_i$ from the collection, 
    there exists $x_j,\; j \neq i,$ such that $l_j$ does not hit 
    $\conv(\{ x_m,\; m \neq i,j \}).$ However, by the assumed properties,
    the only possible choices for such $j$ are $j=i+\varepsilon$ for
    $\varepsilon=+1$ or $\varepsilon=-1$ and the line $l_j$ has to cross
    the segment $\overline{x_i x_{i-\varepsilon}}.$  
    Consequently, each vertex $x_i$ is {\it cut off} from the convex hull
    $\conv(\{x_1,\ldots,x_k\})$ by a line passing through its neighbouring
    vertex and crossing the opposite neighbouring edge. The choices of
    $\varepsilon$ for different vertices are not independent, because if we
    choose to  use $l_i$ to cut off $x_{i+1}$ then $x_{i-1}$ has to be cut off
    by $l_{i-2}.$ Thus, if the number $k$ of
    vertices is odd, this is easily checked to imply that either for all $i=1,\ldots,k$
    the line $l_i$ passes through the segment $\overline{x_i x_{i+1}}$ or for all
    $i=1,\ldots,k$ the line $l_i$ passes through the segment $\overline{x_i x_{i-1}}.$
    In both cases this generates a cycle of order $k$ in ${\cal G}^+[l_1,x_1;\ldots;
    l_k,x_k].$ On the other hand, if $k$ is even, apart from the above two options
    generating a cycle of order $k$ we have another possibility where the collection
    of vertices splits into $k/2$ pairs of neighbours $x_i, x_{i+1}$ such that
    $l_i$ crosses $\overline{x_{i+1} x_{i+2}}$ and $l_{i+1}$ crosses 
    $\overline{x_{i-1} x_i}.$ In this case, however, each such pair of neighbours
    is contained in the trap region generated by any other one and hence we
    get a cycle of order $4$ consisting of some $x_i, y_{i,i+1}, x_j, y_{j,j+1},$
    where  $y_{i,i+1}$ is the intersection point of $l_i$ and $l_{i+1}$ and
    likewise for $y_{j,j+1}.$ In either situation ${\cal G}^+[l_1,x_1;\ldots;l_k,x_k]$
    contains a cycle as required. This completes the proof of Theorem \ref{CORRH}. $\Box$\\ 

   \subsection{Full characterisation of collective factorisation for rectangular fields}\label{FUCH}
    As we already have mentioned above, in the particular case of rectangular Markov fields
    our knowledge is more complete than in the general setting. In fact, we are able to show
    that the necessary precedence graph acyclity condition for collective factorisation
    on all scales, as stated in combined Lemmas \ref{Rysowanie} and \ref{Bezcyklowosc},
    becomes in this context a sufficient condition as well. This is made precise in the
    following theorem.

    \begin{theorem}\label{RectangFactor} 
          Assume that the field ${\cal A}^{\cal M}$ is
          rectangular and let $(l_i,x_i)_{i=1}^k,\; x_i \in l_i,$
          be a collection of pairwise distinct lines and points in general position.
      Then the following are equivalent:
      \begin{enumerate}
       \item The precedence graph  ${\cal G}[l_1,x_1;\ldots;l_k,x_k]$ is acyclic.
       \item The collection $(l_i,x_i)$ admits collective factorisation of correlations
             on all scales.
      \end{enumerate} 
  \end{theorem}
       
  \paragraph{Proof} The implication from 2. to 1. follows directly from Lemmas \ref{Rysowanie}
   and \ref{Bezcyklowosc}. Thus, it only remains to establish the implication from 1. to 2.
   To this end, we make first the following observation, which is also of its own intrinsic
   interest in the context of rectangular fields. 

   Given two non-parallel lines $l$ and $l'$
   on whose parallel translates the activity measure ${\cal M}$ is concentrated, we say that
   a set $D \subseteq {\Bbb R}^2$ is $\{l,l'\}$-convex if for any two points $x,y \in D$
   lying on a common translate of $l$ or $l'$ the entire segment $\overline{xy}$ is contained
   in $D.$ Clearly, $\{l,l'\}$-convexity is a much weaker concept than the usual convexity, in
   particular a $\{l,l'\}$-convex set does not even have to be connected. We shall write
   $\conv_{\{l,l'\}}(A)$ for the $\{l,l'\}$-convex hull of $A \subseteq {\Bbb R}^2,$ that 
   is to say the smallest $\{l,l'\}$-convex set containing $A.$ Now, the crucial observation
   is that in the context of the generalised dynamic representation in Section \ref{GDR}
   specialised for the considered rectangular field, the convexity assumption imposed on the growing
   family $(D_t)$ there can be relaxed to the $\{l,l'\}$-convexity here without any further
   modifications of the theory. Moreover, instead of having the initial set $D_0$ consist
   of a single point, we can now let it consist of any positive finite number of points under
   the $\{l,l'\}$-convexity requirement. Indeed, the proof of Theorem \ref{GEthm} carries
   over verbatim under these relaxed conditions on $(D_t)$ and the $\{l,l'\}$-convexity
   of $(D_t)$ ensures that the intersections of the field lines with $(D_t)$  
   are always connected. It is essential to note at this point that the convexity
   requirement {\it cannot} be lifted for the target domain $D$ though. The point
   is that the consistency statement in Theorem \ref{WlasnosciZgd} does essentially
   require convexity of the domain of the field in its proof and, consequently,
   a rectangular field constructed in some non-convex but $\{l,l'\}$-convex set 
   may fail to coincide with the corresponding restriction of the whole-plane
   consistent field. This phenomenon is closely related to the rather unintuitive
   fact that there exist pairs of compacts $D_1 \subseteq D_2 \subseteq {\Bbb R}^2$
   consisting of a $\{l,l'\}$-convex set $D_1$ and a convex set $D_2$ with the
   property that $D_1$ cannot be extended to $D_2$ through a Hausdorff-continuous
   growing family of $\{l,l'\}$-convex compacts, in sharp contrast to the case
   of convex sets where such an extension is always possible. These issues
   are not discussed in further detail here as falling beyond the scope of
   the present article, yet they are a subject of our active research in
   progress because we believe that the knowledge of the geometry of such
   {\it non-extensible} sets and their corresponding {\it defective} fields
   may shed some further light upon the higher-order correlation structure
   of ${\cal A}^{\cal M}.$           

   To proceed, recall that the vertex set of the precedence graph ${\cal G}[l_1,x_1;\ldots;l_k,x_k]$
   consists of the generating points $x_i$'s and of $y_{i,j}$'s arising as the intersection
   points of the corresponding $l_i$ and $l_j.$
   Since the precedence graph is acyclic, there exists a complete ordering of its vertices
   compatible with the {\it structural} partial order induced by directions of its edges.
   In addition, in view of the minimality of $x_i$'s in the structural order, we can assume that
   no $y_{\cdot,\cdot}$-vertex (intersection vertex) precedes an $x_{\cdot}$-vertex (generating vertex)
   in the considered complete ordering, that is to say the ordered collection of vertices can be written
   as $(x_1,\ldots,x_k;v_1,\ldots,v_m)$ where $v_1,\ldots,v_m$ are all $y_{\cdot,\cdot}$-vertices
   of ${\cal G}[l_1,x_1;\ldots;l_k,x_k].$
   Say that a vertex $u$ of the precedence graph is a direct predecessor of another vertex $w$
   if there is a
   directed edge from $u$ to $w$ containing no other vertices of the graph, which amounts
   to direct precedence of $u$ over $w$ in the structural order on the vertices of the
   precedence graph. Note that
   $y_{\cdot,\cdot}$-vertices have precisely two direct predecessors each, whereas
   $x_{\cdot}$-vertices have no predecessors. We construct inductively a growing family
   $(\tilde{D}_t)$ of $\{l,l'\}$-convex sets in ${\Bbb R}^2$ by putting
   $ \tilde{D}_t := \conv_{\{l,l'\}} D^{\circ}_t,$ where
   \begin{itemize}
    \item We put 
          $ D^{\circ}_0 = \{ x_1,\ldots, x_k\}. $
    \item For $t = i/m,\; 1 \leq i \leq m$ we set 
          $ D^{\circ}_t = \{ x_1,\ldots, x_k, v_1,\ldots,v_i \}.$
    \item For $t \in ((i-1)/m,i/m),\; 1 \leq i \leq m$ we let $w_1,w_2$ be the direct
          predecessors of the vertex $v_i$ (note that $w_1,w_2 \in D^{\circ}_t$ for
          $t \leq (i-1)/m$ by the construction) and we define $D^{\circ}_t$ as the union of
          \begin{itemize}
           \item $D^{\circ}_{(i-1)/m},$
           \item the point $w_1 + m(t-(i-1)/m)(v_{i}-w_1),$
           \item the point $w_2 + m(t-(i-1)/m)(v_{i}-w_2).$
          \end{itemize}
   \end{itemize}       
   In other words, initially the set $\tilde{D}_0$ consists of $\{x_1,\ldots,x_k\}.$
   Since the collection $(l_i,x_i)$ is in general position, we clearly 
   have $\tilde{D}_0 = \conv_{\{l,l'\}} D^{\circ}_0 = 
         \conv_{\{l,l'\}}\{ x_1,\ldots, x_k \} = \{ x_1,\ldots, x_k \}.$
   Next, the intersection points $v_i$ are subsequently added to the
   set $D^{\circ}_t,$ in time intervals of length $1/m.$ With the intersection
   points added, the $\{l,l'\}$-convex hull of the increasing vertex collection
   does no more coincide with the collection itself. The time between the 
   moments of $y_{\cdot,\cdot}$-vertex additions is therefore used to interpolate
   the family $(D^{\circ}_t),$ thus keeping the growth of its $\{l,l'\}$-convex
   hull $(\tilde{D}_t)$ Hausdorff continuous.
   Indeed, the growth of $(D^{\circ}_t)$ is Hausdorff continuous by definition and
   it is easily seen that each new point in $\tilde{D}_{t+dt} \setminus \tilde{D}_t$
   arises as a shift of a point in $D^{\circ}_{t+dt} \setminus D^{\circ}_t$ along
   a translate of either $l$ or $l'$ and, consequently, $\tilde{D}_{t+dt} \setminus \tilde{D}_t$
   is a union of two rectangles of infinitesimal width $O(dt)$ built on two
   possibly degenerated segments parallel to $l$ and $l'$ respectively.
   Note that in general this would not be the case if $v_i$'s were not ordered
   compatibly with the structural ordering \--- adding a single new intersection point
   to $D^{\circ}_t$ might then possibly result in some further intersection
   points, not yet present in $D^{\circ}_t,$ falling into its $\{l,l'\}$-convex
   hull $\tilde{D}_t,$ which might in its turn produce an extra non-degenarate 
   rectangle in $\tilde{D}_t$ thus violating the Hausdorff continuity.
   The importance of imposing on $v_i$'s an ordering compatible with the structural
   order for $(l_i,x_i)_{i=1}^k$ is that it ensures that the intersection points $v_i$
   in $\tilde{D}_t$ are precisely those present in $D^{\circ}_t$ with no extras,
   that is to say for $t \in [i/m,(i+1)/m) \cap [0,1],\; i =0,1,\ldots,m,$ we have
   $$
    \{x_1,\ldots,x_k,v_1,\ldots,v_m\} \cap \tilde{D}_t = 
    \{x_1,\ldots,x_k,v_1,\ldots,v_i\} \cap D^{\circ}_t =
    \{x_1,\ldots,x_k,v_1,\ldots,v_i\}.
   $$  
   In particular, the time order in which the points of the precedence graph
   ${\cal G}[l_1,x_1;\ldots;l_k,x_k]$ show up in $\tilde{D}_t$ is compatible
   with the structural ordering for the collection
   $(l_i,x_i)_{i=1}^k$.
   Observe as well that $\tilde{D}_1$ is a convex set (a rectangle in fact) containing
   all $x_i$'s for $i=1,\ldots,k.$ Moreover, even though $(\tilde{D}_t)$ does not
   satisfy the condition {\bf (D4)}
   in Section \ref{GDR}, it can be easily modified by local smoothing to yield 
   an increasing family $(D_t)$ of $\{l,l'\}$-convex sets satisfying the conditions
   {\bf (D1),(D3),(D4),(D5)} in Section \ref{GDR} and enjoying the property
   (inherited from $(\tilde{D}_t)$) that ${\Bbb A}(l_i) = x_i$ where ${\Bbb A}(\cdot)$
   is the anchor mapping induced by $D_t.$ Furthermore, by the same local smoothing
   modification $D_1$ can be taken convex and containing all $x_i$'s in their interior.
   In view of the discussion on the particular form of the generalised graphical
   construction for rectangular fields and taking into account the {\bf (GE:LineBirth)}
   rule of this construction we conclude the required collective factorisation of
   correlations on all scales, thus completing the proof of Theorem \ref{RectangFactor}. $\Box$

   It should be emphasised at this point that the above argument apparently cannot be repeated
   for more general fields. Indeed, with all directions allowed the condition {\bf (D2)} cannot
   be lifted and $D_0$ has to be a singleton which makes our construction above break down.
   Moreover, in case when only a finite number of, but more than two, directions are allowed,
   the  construction breaks down as well \--- even though a discussion parallel to the above
   can be provided based on the concept of convexity in the directions of the field,
   which does again allow to relax the offending condition {\bf (D2)}, we may in
   general loose the Hausdorff continuity of $\tilde{D}_t$ at the moments when two
   parts of $\tilde{D}_t$ growing in two different directions become for the first
   time connectable along some other direction of the field (never to happen when
   there are only two directions allowed). 


 \subsection{Integral correlation measures}
  With the knowledge of the correlation structure of the field ${\cal A}^{\cal M}$
  as provided so far in this section, we are now in a position to determine
  natural integral correlation measures for the field. To this end, for a bounded
  open convex $D \subseteq {\Bbb R}^2$ and for a planar graph in $D$
  (a collection of edges, not necessarily in $\Gamma_D$) we construct its
  first order directional measure 
  \begin{equation}\label{D1}
   {\cal D}^{(1)}[\gamma] := \sum_{e \in \Edges(\gamma)} \lgth(e) \delta_{l[e]} 
 \end{equation}
  with $\lgth(\cdot)$ denoting the usual Euclidean length and with $\delta_{l[e]}$
  standing for the unit mass at $l[e] \in [[D]].$ The so-defined ${\cal D}^{(1)}(\gamma)$
  is a finite purely atomic measure on $[[D]].$ Likewise, we consider the second order
  directional measure on $[[D]] \times [[D]]$ given by 
  \begin{equation}\label{D2}
   {\cal D}^{(2)}[\gamma] := \sum_{e_1 \in \Edges(\gamma)} \sum_{e_2 \in \Edges(\gamma)}
   \lgth(e_1) \lgth(e_2) \delta_{(l[e_1],l[e_2])}.
  \end{equation} 
  We shall use the expectations of these measures with $\gamma$ drawn from ${\cal A}^{\cal M}_D$
  as natural characteristics of the directional nature of the field. To this end, we put
  \begin{equation}\label{DD1}
   \Delta^{\cal M}_D := {\Bbb E} {\cal D}^{(1)}[{\cal A}^{\cal M}_D].
  \end{equation}
  Write also
  \begin{equation}\label{DD2}
   \Sigma^{\cal M}_D := {\Bbb E} {\cal D}^{(2)}[{\cal A}^{\cal M}_D]
   - \Delta^{\cal M} \otimes \Delta^{\cal M}
  \end{equation}
  where $\otimes$ denotes the measure producting operation. Note that $\Sigma^{\cal M}_D$
  has its natural interpretation as the directional covariance measure. Below we explicitly
  determine both the first and second order directional measures for ${\cal A}^{\cal M}.$ 

  \begin{corollary}\label{FOD}
   For $l \in [[D]]$ we have
   $$ \Delta^{\cal M}_D(dl) = \lgth(l \cap D) {\cal M}(dl). $$
  \end{corollary}

  \paragraph{Proof}
   For each $l \in [[D]]$ we have
   $$ \Delta^{\cal M}(dl) = \int_{l \cap D} \sigma^{\cal M}[dl,x] \ell(dx) $$
   where $\ell$ is the length element. Our statement follows now by Theorem \ref{Corr12}. $\Box$
 
  \begin{corollary}\label{SOD}
   For $l_1,l_2 \in [[D]]$ we have
   $$ \Sigma^{\cal M}_D(dl_1,dl_2) = \left\{ \begin{array}{ll}
      \int_{l \cap D} \int_{l \cap D} \exp(-2 {\cal M}([[\overline{xy}]]))
      \lgth(dx) \lgth(dy) {\cal M}(dl), & \mbox{ if } l_1=l_2=l,\\
      0, & \mbox{ otherwise.} \end{array} \right. $$
  \end{corollary}  

  \paragraph{Proof}
   For $l_1, l_2 \in [[D]]$ we have
   $$ \Sigma^{\cal M}_D(dl_1,dl_2) = \int_{l_1 \cap D} \int_{l_2 \cap D}
      \sigma^{\cal M}[dl_1,x_1;dl_2,x_2] \ell(dx_2) \ell(dx_1) - 
      \Delta^{\cal M}_D(dl_1) \Delta^{\cal M}_D(dl_2). $$ Our assertion
   follows now by Theorem \ref{Corr12}. $\Box$

\subsection{Alternative correlation measures}
  In addition to the events that there exist field lines passing through given
  points in given directions as considered in the definition of $\sigma^{\cal M}[\ldots]$
  above, one can also study events that there exist field angles with vertices at
  given points and with their arms extending along given half-lines.
  In formal terms, we put
  \begin{equation}\label{CORRA}
  \varsigma^{\cal M}[d\angle_1;\ldots;d\angle_k] := 
   d{\Bbb P}\left( \forall_{i=1}^k \angle_k \in \Angles({\cal A}^{\cal M})\right),
  \end{equation}
  where $\Angles(\gamma),\; \gamma \in \Gamma,$ denotes the collection
  of convex angles in the admissible polygonal graph $\gamma$ and where 
  $\angle_i := \angle[\vec{l}_i,\vec{l'}_i],\; i=1,\ldots,k,$ stands
  for the convex angle between two directed half-lines $\vec{l}_i$
  and $\vec{l'}_i$ outgoing from the vertex of the angle and directed
  away from the vertex. It should be noted at this point for formal
  completeness that the concave angles arising as the complements of
  the convex angles of  $\gamma$ {\it are not} included into
  $\Angles(\gamma).$ In analogy to Lemma \ref{CORRHC}, a similar 
  factorisation result for the probabilities of such events can be concluded
  from the generalised dynamic representation.
  \begin{corollary}\label{CORRAN}
   Assume the collection $\angle_1,\ldots,\angle_k$ is such that the vertices
   $x_i$ of $\angle_i$ are all different and such that the interior of 
   the convex angle $\angle_1$ is disjoint with the convex hull 
   $\conv(\{x_2,\ldots,x_k\}).$ Then
   $$ \varsigma^{\cal M}[d\angle_1;\ldots;d\angle_k] = 
      {\cal M}(dl_1){\cal M}(dl'_1) \varsigma^{\cal M}[d\angle_2;\ldots;d\angle_k]. $$
  \end{corollary}
  \paragraph{Proof}
   Indeed, under the assumptions of the corollary the increasing family $D_t$ in
   the generalised dynamic representation in Section \ref{GDR} can be chosen so that it
   first hits all $x_2,\ldots,x_k$ and finally it reaches $x_1$ before hitting any other
   point of $\angle_1.$ Our claim follows now by the rule {\bf (GE:VertexBirth)}
   of the generalised representation. $\Box$          
   
  Suppose now that we assign alternating labels $+1$ and $-1$ to regions separated
  by the contours of the field which thus become interfaces between $+1$ and $-1$
  phases. This gives rise to two possible label assignments in ${\Bbb R}^2$ and
  we pick any of them with probability $1/2.$ Write
  $[{\cal A}^{\cal M}]_x$ for the label assigned by ${\cal A}^{\cal M}$ 
  at $x \in {\Bbb R}^2.$ Then for two points $x,y \in {\Bbb R}^2$ we can define
  the two-point label-correlation function
  \begin{equation}\label{PPcorr}
   \rho^{\cal M}_{x,y} := {\Bbb E}  [{\cal A}^{\cal M}]_x
    [{\cal A}^{\cal M}]_y  - {\Bbb E} [{\cal A}^{\cal M}]_x 
    {\Bbb E} [{\cal A}^{\cal M}]_y.
  \end{equation}
  It is easily seen that, writing $n_{x,y} := \card (\overline{xy} \cap {\cal A}^{\cal M})$ 
  for the number of edges of the field ${\cal A}^{\cal M}$ crossing the segment $\overline{xy},$ 
  $$ \rho^{\cal M}_{x,y} = {\Bbb P}(n_{x,y} \mbox{ is even }) -
     {\Bbb P}(n_{x,y} \mbox{ is odd }). $$
  Recalling from the {\bf Linear sections} statement of Theorem \ref{WlasnosciZgd}
  that $n_{x,y}$ is Poisson with parameter ${\cal M}([[\overline{xy}]])$
  with $\overline{xy}$ denoting the segment joining $x$ to $y,$ we conclude that
  \begin{corollary}\label{PPcorrC}
   We have
   $$ \rho^{\cal M}_{x,y} = \exp(-{\cal M}([[\overline{xy}]])). $$
  \end{corollary}
  Roughly speaking, this shows that for regularly behaved activity measures
  ${\cal M}$ the two-point label-correlation functions of the field
  ${\cal A}^{\cal M}$  exhibit exponential decay with the distance between
  $x$ and $y$ which places the inverse temperature $\beta = 1$ corresponding
  to the consistent regime in the high temperature regime for polygonal fields,
  above the phase transition point.


\section{Disagreement loop death and birth dynamics}\label{DLOOP}
 In this section we discuss a random dynamics on the space $\Gamma_D$ of admissible
 polygonal configurations which leaves the law of the consistent field ${\cal A}^{\cal M}_D$
 invariant, with $D$ standing as usual for an open bounded convex set in ${\Bbb R}^2.$
 We also provide a modification of this dynamics which can be used for Monte-Carlo
 simulation of ${\cal A}^{{\cal M};\beta}_D$ for all $\beta \in {\Bbb R}.$ We build upon
 \cite{SC05} in our presentation of the standard dynamics based on an important
 concept of a {\it disagreement loop}.

 \subsection{Standard disagreement loop dynamics for consistent regime}
 To proceed we place ourselves within the context of the standard dynamic
 representation as given in Subsection \ref{DYNRES} and suppose that we
 observe a particular realisation ${\gamma} \in {\Gamma}_D$ of the
 polygonal field ${\cal A}^{\cal M}_D$ and that we modify the
 configuration by adding an extra birth site $x_0$ to the existing
 collection of birth sites for $\gamma,$ while keeping the evolution rules
 {\bf (E1-4)} for all the particles, including
 the two newly added ones if $x_0 \in D$ and the single newly
 added one if $x_0 \in \partial D.$ Denote the resulting new (random)
 polygonal configuration by $\gamma \oplus x_0.$ A simple
 yet crucial observation is that for $x_0 \in D$ the symmetric difference
 $\gamma \triangle [\gamma \oplus x_0]$ is almost surely a single loop
 (a closed polygonal curve), possibly self-intersecting and possibly
 chopped off by the boundary. Indeed, this is seen as follows.
 The leftmost point of the loop $\gamma \triangle [\gamma \oplus x_0]$
 is of course $x_0.$ Each of the two {\it new} particles $p_1, p_2$
 emitted from $x_0$ move independently, according to ${\bf (E1-4)},$ each
 giving  rise to a {\it disagreement path}. The initial segments of such a
 disagreement path correspond to the movement of a particle, say  $p_1,$
 before its annihilation in the first collision. If this is a collision
 with the boundary, the disagreement path gets chopped off and terminates
 there. If this is a collision with a segment of the original configuration
 $\gamma$ corresponding to a certain {\it old} particle $p_3,$ the
 {\it new} particle $p_1$ dies but the disagreement path continues
 along the part of the trajectory of $p_3$ which is contained in
 $\gamma$ but not in $\gamma \oplus x_0.$ At some further moment $p_3$
 dies itself in $\gamma,$ touching the boundary or killing another
 particle $p_4$ in $\gamma.$ In the second case, however, this collision
 only happens for $\gamma$ and not for $\gamma \oplus x_0$ so the
 particle $p_4$ survives (for some time) in $\gamma \oplus x_0$
 yielding a further connected portion of the disagreement path initiated by $p_1,$
 which is contained in $\gamma \oplus x_0$ but not in $\gamma$ etc.
 A recursive continuation of this construction shows that the disagreement
 path initiated by $p_1$ consists alternately of connected polygonal subpaths
 contained in $[\gamma \oplus x_0] \setminus \gamma$ (call these {\it creation phase}
 subpaths) and in $\gamma \setminus [\gamma \oplus x_0]$ (call these {\it annihilation
 phase} subpaths). Note that this disagreement path is self-avoiding and, in
 fact, it can be represented as the graph of some piecewise linear
 function $t \mapsto y(t).$ Clearly, the same applies for the disagreement path
 initiated by $p_2.$ An important observation is that whenever two {\it creation phase}
 or two {\it annihilation phase} subpaths of the two disagreement paths hit each other,
 both disagreement paths die at this point and the disagreement loop closes
 (as opposed to intersections of segments of different phases which do not
 have this effect). Obviously, if the disagreement loop does not close
 in the above way, it gets eventually chopped off by the boundary.
 We shall write $\Delta^{\oplus}[x_0;\gamma] = \gamma \triangle [\gamma \oplus x_0]$
 to denote the (random) disagreement loop constructed above. It remains
 to consider the case $x_0 \in \partial D,$ which is much simpler
 because there is only one particle emitted and so $\Delta^{\oplus}[x_0;\gamma]
 = \gamma \triangle [\gamma \oplus x_0]$ is a single self-avoiding
 polygonal path eventually chopped off by the boundary. We shall often abuse
 the language calling such a disagreement path $\Delta^{\oplus}[x_0;\gamma]$
 a (degenerate) disagreement loop as well.

 Likewise, a disagreement loop (or path) arises if we {\it remove} one birth site
 $x_0$ from the collection of birth sites of an admissible polygonal configuration
 $\gamma \in \Gamma_D,$ while keeping the evolution rules for
 all the remaining particles. We write $\gamma \ominus x_0$ for the
 configuration obtained from $\gamma$ by removing $x_0$ from the list
 of the birth sites, while the resulting random disagreement loop (or path) is
 denoted by $\Delta^{\ominus}[x_0;\gamma]$ so that $\Delta^{\ominus}[x_0;\gamma]
 = \gamma \triangle [\gamma \ominus x_0].$ We refer the reader to Section
 2.1 in \cite{SC05} for further discussion.

 With the above terminology we are in a position to describe a random
 dynamics on the configuration space $\Gamma_D,$
 which leaves invariant the law of the
 consistent polygonal process ${\cal A}^{\cal M}_D.$ Particular care
 is needed, however, to distinguish between the notion of time considered in
 the dynamic representation of the field as well as throughout the construction
 of the disagreement loops above, and the notion of time to be introduced
 for the random dynamics on $\Gamma_D$ constructed below. To make
 this distinction clear we shall refer to the former as to the {\it
 representation time} (r-time for short) and shall reserve for it the notation
 $t,$ while the latter will be called the {\it simulation time} (s-time
 for short) and will be consequently denoted by $s$ in the sequel.

 Consider the following pure jump birth and death type Markovian
 dynamics on $\Gamma_D,$ with $\gamma_s = \gamma_s^D$ standing
 for the current configuration
 \begin{description}
  \item{\bf DL:birth} With intensity 
   $\lnaw {\cal M} \pnaw(dx) ds$ for $x \in D$ and 
   with intensity $\lnaw {\cal M};\partial D \pnaw(dx) ds$ for $x \in \partial D$ 
   set $\gamma_{s+ds} := \gamma_s \oplus x,$
  \item{\bf DL:death} For each birth site $x$ in $\gamma_s$
   with intensity $ds$ set $\gamma_{s+ds} := \gamma_s \ominus x.$
 \end{description}
 If none of the above updates occurs we keep $\gamma_{s+ds} =
 \gamma_s.$ It is convenient to perceive the above dynamics
 in terms of generating random disagreement loops $\lambda$ and
 setting $\gamma_{s+ds} := \gamma_s \triangle \lambda,$ with the
 loops of the type $\Delta^{\oplus}[\cdot,\cdot]$ corresponding
 to the rule {\bf DL:birth} and $\Delta^{\ominus}[\cdot,\cdot]$
 to the rule {\bf DL:death}.

 As a direct consequence of the dynamic representation
 of the consistent field ${\cal A}^{\cal M}_D$ we obtain
 \begin{theorem}\label{AR1}
  The distribution of the polygonal field
  ${\cal A}^{\cal M}_D$ is the unique invariant law of the
  dynamics given by {\bf DL:birth} and {\bf DL:death}. The resulting
  s-time stationary process is reversible. Moreover, for any initial
  distribution of $\gamma_0$ the laws of the polygonal fields $\gamma_s$
  converge in variational distance to the law of ${\cal A}^{\cal M}_D$
  as $s \to \infty.$
 \end{theorem}
 The uniqueness and convergence statements in the above theorem
 require a short justification. They both follow by the observation
 that, in finite volume, regardless of the initial state, the process
 $\gamma_s$ spends a non-null fraction of time in the empty state
 (no polygonal contours).
 Indeed, this observation allows us to conclude the required uniqueness
 and convergence by a standard coupling argument, e.g. along the lines 
 of the proof of Theorem 1.2 in \cite{LG}. 

\subsection{Standard disagreement loop dynamics for general temperatures}\label{GenTemp}
 Below, we  show that the laws of the Gibbs-modified polygonal fields
 ${\cal A}^{{\cal M};\beta}_D,\; \beta \in {\Bbb R},$ arise as the
 unique invariant distributions for appropriate modifications of the
 reference dynamics {\bf DL:birth, DL:death}. The main change is
 that the birth and death updates are no more performed unconditionally,
 they pass an {\it acceptance test} instead and are accepted
 with certain state-dependent probabilities whereas upon failure of
 the acceptance test the update is discarded. Consider the
 following dynamics
 \begin{description}
  \item{${\bf DL:birth[\beta]}$}
   With intensity 
   $\lnaw {\cal M} \pnaw(dx) ds$ for $x \in D$ and 
   with intensity 
   $\lnaw {\cal M};\partial D \pnaw(dx) ds$ for $x \in \partial D$ 
   do
   \begin{itemize}
    \item put $\delta := \gamma_s \oplus x,$
    \item accept $\delta$ with probability 
      $\exp(-[\beta-1] (L^{\cal M}(\delta) - L^{\cal M}(\gamma_s)))$
      if $L^{\cal M}(\delta) > L^{\cal M}(\gamma_s)$ and with probability
      1 otherwise,
    \item if accepted, set $\gamma_{s+ds} := \delta,$ otherwise
          keep $\gamma_{s+ds} := \gamma_s.$
   \end{itemize}
  \item{${\bf DL:death[\beta]}$}
   For each birth site $x$ in $\gamma_s$ with intensity $ds$ do
   \begin{itemize}
    \item put $\delta := \gamma_s \ominus x,$
    \item accept $\delta$ with probability 
      $\exp(-[\beta-1] (L^{\cal M}(\delta) - L^{\cal M}(\gamma_s)))$
      if $L^{\cal M}(\delta) > L^{\cal M}(\gamma_s)$ and with probability
      1 otherwise,
    \item if accepted, set $\gamma_{s+ds} := \delta,$ otherwise
          keep $\gamma_{s+ds} := \gamma_s.$
   \end{itemize}
 \end{description}
 In analogy with its original reference form {\bf DL:birth, DL:death},
 the above dynamics should be thought of as generating random disagreement
 loops $\lambda$ and setting $\gamma_{s+ds} := \gamma \triangle \lambda$
 provided $\lambda$ passes the acceptance test.
 The following theorem justifies the above construction.
 \begin{theorem}\label{SYMULACJADL}
  For each $\beta \in {\Bbb R}$ the law of the polygonal process 
  ${\cal A}^{{\cal M};\beta}_D$ is the unique invariant
  distribution of the dynamics ${\bf DL:birth[\beta],
  DL:death[\beta]}.$ The resulting s-time stationary process is reversible.
  For any initial distribution of $\gamma_0$ the laws of the polygonal
  fields $\gamma_s$ converge in variational distance to the law of
  ${\cal A}^{{\cal M};\beta}_D$ as $s \to \infty.$
 \end{theorem}
 Indeed, the invariance follows by a straightforward check of usual
 detailed balance conditions whereas the convergence and uniqueness
 statements are verified as in Theorem \ref{AR1}.  

\subsection{Generalised disagreement loop dynamics}\label{GDLD}
 Even though the above disagreement loop dynamics in its standard version
 has been constructed using the standard dynamic representation of
 Subsection \ref{DYNRES}, exactly the same can be made for the generalised
 dynamic  representation of Section \ref{GDR}. As easily verified, this also
 leads to disagreement  loop and path creation and annihilation, the only
 difference being that boundary birth events in the standard set-up are replaced
 by general line birth events in context of the generalised dynamic construction.
 Vertex creation and annihilation yield disagreement loops whereas line creation
 and annihilation yield disagreement paths.
 Thus, the {\bf DL:birth} and {\bf DL:death} moves get replaced by
 \begin{description}
  \item{\bf GenDL:birth} With intensity $\lnaw {\cal M} \pnaw(dx) ds$ for $x \in D$
    set $\gamma_{s+ds} := \gamma_s \oplus x.$ With intensity ${\cal M}(dl)ds$ for
    $l \in [[D]]$ set $\gamma_{s+ds} := \gamma_s \oplus l,$ where $\gamma_s \oplus
    l$ arises from $\gamma_s$ by letting $l$ be born at its anchor point ${\Bbb A}(l)$
    and making it  extend thereupon according to the evolution rules of the generalised
    construction, which gives rise to a disagreement path.
  \item{\bf GenDL:death} For each point birth site $x$ in $\gamma_s$
   with intensity $ds$ set $\gamma_{s+ds} := \gamma_s \ominus x.$
   For each line birth site ${\Bbb A}(l)$ of a line $l$ in $\gamma_s,$
   with intensity $ds$ set $\gamma_{s+ds} := \gamma_s \ominus l,$
   where again $\gamma_s \ominus l$ arises from $\gamma_s$ by killing
   at ${\Bbb A}(l)$ the field line $l$ of $\gamma_s,$ which gives rise
   to a disagreement path.    
 \end{description}
 In full analogy with Theorem \ref{AR1} we have 
 \begin{theorem}\label{GenAR1}
  The distribution of the polygonal field
  ${\cal A}^{\cal M}_D$ is the unique invariant law of the
  dynamics given by {\bf GenDL:birth} and {\bf GenDL:death}. The resulting
  s-time stationary process is reversible. Moreover, for any initial
  distribution of $\gamma_0$ the laws of the polygonal fields $\gamma_s$
  converge in variational distance to the law of ${\cal A}^{\cal M}_D$
  as $s \to \infty.$
 \end{theorem}
 The modification ${\bf GenDL[\beta]}$ of the {\bf GenDL} dynamics for general
 $\beta \in {\Bbb R}$ goes along exactly the same lines as the construction of
 ${\bf DL[\beta]}$
 from {\bf DL}, by introducing suitable acceptance tests,
 and hence we omit the standard details here. Obviously, the 
 corresponding version of Theorem \ref{SYMULACJADL} holds in full analogy.
 Clearly, the so-defined dynamics depends then on the choice of the
 increasing family $D_t.$ Note that hybrid dynamics can also be
 considered, switching between different choices of $D_t$ in the course
 of the simulation time $s,$ which is sometimes beneficial e.g. in
 image processing applications, see \cite{SL08}.    


\subsection{Defective disagreement loop dynamics}\label{DDLD}
 To conclude the discussion of disagreement loop dynamics we introduce one more
 variant thereof, which we call {\it defective dynamics}. The idea stems from the
 observation that the standard dynamic construction in Subsection \ref{DYNRES}
 and its generalised version in Section \ref{GDR} rely
 respectively on the time axis direction and on the anchor mapping ${\Bbb A}:[[D]] \to D$
 to decide {\it in which of the two possible directions to proceed along a new line}
 each time a directional update occurs and a new edge is to be created. This decision
 mechanism is then inherited by the disagreement loop dynamics where the current
 direction of an unfolding disagreement loop changes in the following situations:
 \begin{enumerate}
  \item upon a directional update during a creation phase,
  \item when a directional update is present on the path being erased during
        an annihilation phase,
  \item when switching between creation and annihilation phases.
 \end{enumerate}
 In all these situations the choice of the direction to follow along the
 new line is always the same: forward in time in the standard case and 
 away from the corresponding anchor point in the generalised case. 
 In certain contexts there may be good reasons to try to overcome these
 restrictions and to allow for more general directional decision
 strategies even at the cost of obtaining disagreement loop dynamics
 which do not correspond to any well-defined dynamic representations.
 One of such situations arises when studying the correlation
 functions of the polygonal field as defined in Section \ref{CORCR} and we
 shall indeed use the defective construction to provide below a proof of
 exponential correlation decay for rectangular fields, see Section
 \ref{CORRDEC}. Another setting where more general directional
 decision strategies may be of use are image processing applications
 of polygonal fields, see \cite{SL08}.

 The idea underlying the defective disagreement loop dynamics is hence to keep the
 crucial feature of the above {\bf DL} and {\bf GenDL} dynamics, being that
 the symmetric difference between a configuration and its update is always a
 disagreement loop algorithmically constructible in the course of alternating
 creation and annihilation phases, but to give up the requirement that the
 disagreement loop generating mechanism be determined by a dynamic representation.  
 Such an approach makes this {\it disagreement loop generating mechanism} the
 core of the construction. In our construction of the defective dynamics we shall
 assume that upon each directional update along an unfolding disagreement loop, the
 decision on which of the two possible directions to choose along the new line is
 taken as {\it a deterministic function of the current portion of the disagreement
 subpath built by the time of the update}.
 To put it in formal terms, assume we are given a measurable {\it directional
 decision} mapping/rule $\phi = \phi_D$ assigning to each connected non-closed but
 possibly self-intersecting polygonal path $\lambda$ contained in the domain
 $D$ of the field and to a line $l \in [[D]]$ passing through one of the two
 ends of $\lambda,$ one of the two possible directions along $l.$ In our dynamics
 below $\phi(\lambda,l)$ will be used to determine the direction in which the
 disagreement path $\lambda$ is to extend upon a directional update into $l.$
 For a polygonal configuration $\gamma \in \Gamma_D,$ a straight line $l \in [[D]]$
 assumed not to be an extension of an edge of $\gamma$ and a point
 $x \in l \setminus \gamma$ we construct the configuration  $\gamma \oplus (l,x)
 \in \Gamma_D,$ interpreted as the result of adding to $\gamma$ a {\it linear germ}
 $(l,x)$ born at $x$ and initially evolving along $l \ni x$ or, in other words,
 creating a new field line for $\gamma$ at $x.$ The construction goes as
 follows. We add a line birth site for $l$ at $x$ at the time $t = 0$ whereupon
 we let the newly born {\it linear germ} extend in both directions along $l$ away from
 $x,$ say with constant unit speed for definiteness although other options are also
 possible as discussed below. While extending, the resulting new edge is subject to
 directional updates: for each line $l'$ which the growing edge crosses the probability
 that a directional update occurs and a new growth direction along $l'$ is assumed
 is ${\cal M}(dl')$ \--- if this happens at some time moment $t,$ the new direction
 along $l'$ is always chosen according to the directional decision rule
 $\phi(\lambda_t,l')$ and the growth continues with constant unit speed,
 with $\lambda_t$ standing for the current path constructed by the time
 $t.$ At some point, the resulting polygonal path may hit
 the boundary of $D$ or an existing edge $e$ of $\gamma,$ in which case it stops
 growing in this direction. Upon a collision with an edge we switch from {\it creation}
 mode to {\it annihilation} mode and proceed (say, again with constant unit speed)
 from the collision point along the edge $e$
 in the direction indicated by $\phi(\lambda_t,l[e]),$ erasing this part of the edge.
 This continues till we reach the boundary $\partial D$ or another edge $e'$ of
 $\gamma.$ In the latter case, if the edge $e'$ extends in the direction agreeing
 with $\phi(\lambda_t,l[e'])$ from the collision point, we proceed along $e'$
 and keep erasing it.
 If, however, the edge $e'$ extends from the collision point in
 the direction opposite to $\phi(\lambda_t,l[e']),$ we switch back to the {\it creation}
 mode and start extending $e'$ in the direction indicated by $\phi(\lambda_t,l[e']).$
 We continue this procedure for both {\it branches} of the linear germ emitted from the
 line birth site $(l,x).$ Both branches give rise to directed polygonal {\it disagreement
 paths} constituting the symmetric difference $\gamma \triangle (\gamma \oplus (l,x)),$
 both consisting of alternating polygonal subpaths corresponding to creation and annihilation
 stages of the construction and respectively referred to as {\it creation phase} and
 {\it annihilation phase} subpaths in the sequel. The procedure terminates when both branches
 get killed in collision with $\partial D.$ It may be useful to observe that the requirement
 that both branches of the disagreement loop grow with constant unit speed has no effect
 on the final shape of individual branches but it determines the precedence of their
 intersection points, should they meet. Thus, in principle, the unit growth speed may
 be alternatively replaced by other natural options such as letting the right (left)
 branch grow first till it reaches the boundary and to let the left (right) branch
 grow thereupon, or even admitting different growth speeds for both branches or speeds
 depending on the current form of the disagreement loop under construction.
 Below we shall always assume one of these options is chosen and fixed,
 under each such choice the discussion below remains valid.
 The co-existence of two disagreement branches obliges us to impose one  more rule on the
 construction of $\gamma \oplus (l,x):$ if a creation phase subpath of one
 branch hits a creation phase subpath of another branch, both branches
 get killed at this point, referred to as the {\it cut-off point} below,
 their continuation possibly already constructed gets chopped off and
 we consider the disagreement loop complete thus terminating the construction;
 likewise for two annihilation phase subpaths of different branches colliding.
 If two subpaths of different phases coming from different branches meet,
 nothing happens though. Note that there can be several such cut-off points
 and the disagreement loop construction terminates at the {\it first} cut-off
 point with respect to the ordering imposed by the time flow \--- consequently
 the choice of a particular growth speed protocol as discussed above may affect
 the choice of the loop-closing cut-off point and hence also the final shape
 of the disagreement loop.    

 So far, this construction looks very similar to the construction
 of the disagreement paths and loops in {\bf DL} and {\bf GenDL} above. There is one crucial
 difference though: it may happen here that {\it a self-intersection occurs along one branch
 of the disagreement path being constructed} \--- this could never have happened in the
 disagreement loop dynamics {\bf DL} and {\bf GenDL} corresponding to dynamic representations
 because of
 the temporal ordering along branches of disagreement paths/loops there, and this is
 the price we pay here for the freedom of choice of the directional decision mapping
 $\phi(\cdot,\cdot).$ If such a self-intersection occurs past the cut-off point, we 
 just ignore it. If it happens before the cut-off, we are in trouble though. 
 A look at the description of the disagreement path construction
 procedure here readily shows that the presence of such a self-intersection
 and hence of a cycle along a single disagreement branch before the cut-off means
 \begin{itemize}
  \item either that at this point we start discarding the changes previously introduced
        by erasing segments previously created (meeting of two creation phases)
  \item or that the $\oplus (l,x)$ operation is not reversible (creation phase edge
        hits annihilation phase subpath preceding it along the disagreement branch)
        in that the disagreement loop $\gamma \triangle (\gamma \oplus (l,x)),$ even
        though it can be used for the transition $\gamma \to (\gamma \oplus (l,x)),$ is no
        more a valid disagreement loop upon switching the creation and annihilation
        phases and thus cannot be used for the reverse transition
        $(\gamma \oplus (l,x)) \to \gamma.$
 \end{itemize}
 The remaining two situations: intersection of two annihilation subpaths or an
 annihilation edge hitting a creation subpath preceding it along the disagreement
 branch, cannot occur in our construction as this would either imply the existence
 of a node of order higher than $2$ in the original configuration $\gamma$
 (in the first case) or stand in contradiction to the disagreement loop
 construction rules (in the second case). 
 To prevent the above pathologies from happening we say that the construction
 of the disagreement loop/path {\it fails} if a cycle occurs along one of the 
 branches. Consequently, the configuration $\gamma \oplus (l,x)$ is
 not always defined, hence the name {\it defective dynamics}.

 In full analogy with the above definition of the $\oplus$ operation, we can define
 the configuration $\gamma \ominus (l,x) \in \Gamma_D$
 for line $l$ extending an edge $e \in \Edges(\gamma),\; \gamma \in \Gamma_D,$
 with $x \in e,$ which is to be interpreted as initiating the annihilation
 phase in $\gamma$ from the linear germ $(l,x)$ or, in other words,
 killing/annihilating the field line $l$ of $\gamma$ at the point $x.$
 This operation also gives rise to a disagreement loop and its outcome may
 be undefined as well, on equal rights with $\ominus.$
 
 The operations $\oplus (l,x)$ and $\ominus (l,x)$ combine the properties
 of line birth/death and vertex birth/death operations of the generalised
 disagreement loop dynamics in Subsection \ref{GDLD} and there are no separate $\oplus l,
 \oplus x,\ominus l,\ominus x$ operations in the defective dynamics.
 

 Directly from the definition of the line creation operation $\oplus$ it follows that if
 for some $\gamma \in \Gamma_D$ and a linear germ $(l,x),\; x \in l \setminus \gamma,$
 another configuration $\gamma' \in \Gamma_D$ is reachable from $\gamma$ upon adding the
 linear germ $(l,x),$ that is to say $\gamma'$ is a possible outcome for
 $\gamma \oplus (l,x),$ then the probability element that $\gamma \oplus (l,x) = \gamma'$ is
 \begin{equation}\label{TransProb1}
   \exp(-L^{\cal M}(\gamma' \setminus \gamma)) 
    \prod_{e \in \Edges(\gamma' \setminus \gamma),\; e \not\ni x} {\cal M}(dl[e]), 
 \end{equation}
 where $\Edges(\gamma' \setminus \gamma)$ stands for the collection of edges
 of $\gamma'$ {\it which are not extensions of edges of} $\gamma$ and thus they
 have to result from directional updates during creation phase. Indeed, the product
 $\prod {\cal M}(dl[e])$ in (\ref{TransProb1}) comes from the probability cost of
 creation phase directional updates  yielding all edges in
 $\Edges(\gamma' \setminus \gamma)$ save the initial one
 germinating from $x,$ whereas $\exp(-L^{\cal M}(\gamma' \setminus \gamma))$ is due to
 the absence of directional updates along the created polygonal paths of
 $\gamma' \setminus \gamma.$
 In a similar way we see that if for some $(l,x)$ the configuration $\gamma \in \Gamma_D$
 contains an edge $e$ with $x \in e$ and $l = l[e]$ then for another configuration $\gamma'$
 reachable from $\gamma$ by annihilating the line $l$ at $x$ we have the probability
 element that $\gamma \ominus (l,x) = \gamma'$ given by 
 \begin{equation}\label{TransProb2}  
    \exp(-L^{\cal M}(\gamma' \setminus \gamma)) 
    \prod_{e \in \Edges(\gamma' \setminus \gamma)} {\cal M}(dl[e]).
 \end{equation}
 It is also crucial to observe that each valid transition
 $\gamma \to \gamma' = \gamma \oplus (l,x)$ admits its inverse
 obtainable as $\gamma' \to \gamma = \gamma' \ominus (l,x)$ under
 the same directional decisions made in the course of the
 disagreement loop construction but with the creation and
 annihilation phases switched \--- this is due to the fact that
 the directional decision rule $\phi(\cdot,\cdot)$ only takes
 into account the current state of the unfolding disagreement
 loop (which is the same for $\gamma \to \gamma'$ and 
 $\gamma' \to \gamma$ as the symmetric difference) rather than
 the actual configuration (which differs for both transitions). 

 Consider now a measurable set $A_L \subseteq [[D]]$ and a measurable mapping
 ${\Bbb A}^*: A_L \to D$ such that ${\Bbb A}^*(l) \in l.$ These are interpreted
 respectively as the set of lines allowed for creation/annihilation and the assignment to
 each line of its germination point. The mapping ${\Bbb A}^*$ can also be
 randomised as long as it is kept independent of the dynamics below, but
 for simplicity we do not discuss the details of this straightforward extension
 as it is not going to be used in this paper. Define the following
 {\it defective disagreement loop dynamics} unfolding in s-time:
 \begin{description}
  \item{\bf DefDL:birth} With intensity ${\cal M}(dl)ds$ for $l \in A_L$
    set $\gamma_{s+ds} := \gamma_s \oplus (l,{\Bbb A}^*(l)).$ Should
    the proposed update fail, keep $\gamma_{s+ds} = \gamma_s.$
  \item{\bf DefDL:death} For each line $l \in A_L$ extending an edge $e$ of 
   $\gamma_s$ such that ${\Bbb A}^*(l) \in e,$ with intensity $ds$
   set $\gamma_{s+ds} := \gamma_s \ominus (l,{\Bbb A}^*(l)).$
   Should the proposed update fail, keep $\gamma_{s+ds} = \gamma_s.$       
 \end{description}  
 Recalling the above transition probabibilities (\ref{TransProb1})
 and (\ref{TransProb2}) corresponding to $\oplus$ and $\ominus$ moves
 and taking into account the discussion following their displays 
 allows for a straightforward verification of the detailed balance
 conditions which readily yields
 \begin{theorem}\label{DefAR1}
  The distribution of the polygonal field ${\cal A}^{\cal M}_D$ is an
  invariant law for the dynamics given by {\bf DefDL:birth} and
  {\bf DefDL:death}. The resulting s-time stationary process is
  reversible.
 \end{theorem}  
 Note that in strong contrast to analogous Theorems \ref{AR1} and
 \ref{GenAR1} we do not state uniqueness or convergence to the law of
 the polygonal field ${\cal A}^{\cal M}_D$ here and this is because
 in general neither of them holds. Indeed, this is due to the fact that
 a defective dynamics does not have to be transitive: there may exists
 pairs of admissible polygonal configurations which cannot be reached
 from each other by $\oplus$ and $\ominus$ moves. This is closely
 related to the possible failures of these updates and it is the price
 to pay for the free choice of the directional decision mapping. The
 structure of the recurrent classes of defective dynamics on $\Gamma_D$
 is a subject of our research in progress. Observe as well that, in
 view of the discussion following the displays (\ref{TransProb1}) and
 (\ref{TransProb2}), the requirement that the directional decision
 rule $\phi(\cdot,\cdot)$ only depends on the current state of the
 disagreement loop rather than on the entire actual configuration
 seems indispensable for the detailed balance and the reversibility
 of the dynamics.  

 \paragraph{Variants and applications of the defective dynamics}
  The modification ${\bf DefDL[\beta]}$ of the defective disagreement loop
  dynamics {\bf DefDL} for general $\beta \neq 1$ is straightforward and
  goes by introducing suitable acceptance test for update proposals 
  in full analogy with the corresponding modification of the ${\bf DL}$
  dynamics aimed at obtaining ${\bf DL[\beta]}$ in Subsection \ref{GenTemp}.  

  The class of possible directional decision rules is extremely rich and
  various strategies can be chosen encompassing different goals. One possible
  option of constructing the directional decision mapping goes by choosing
  an arbitrary measurable {\it generalised anchor mapping} ${\Bbb A}:[[D]] \to D$
  with ${\Bbb A}(l) \in l$ to be interpreted as the {\it initial point} of the
  line $l.$ Then the direction $\phi(\lambda,l)$ the disagreement loop $\lambda$
  is to assume for its unfolding along $l$ upon a directional update is always
  set {\it away from the anchor point} ${\Bbb A}(l).$ Note that if the
  anchor mapping corresponds to an increasing family $(D_t)$ of convex 
  compacts satisfying {\bf (D1-5)} as in Section \ref{GDR} then the 
  resulting defective dynamics $\oplus$ and $\ominus$ operations
  coincide respectively with the line birth and line death operations of
  the generalised diasgreement loop dynamics described in Subsection
  \ref{GDLD}.  

  The directional decision mapping can also be randomised, as long as 
  $\phi(\lambda_t,\cdot)$ is independent of the current field
  configuration conditionally given $\lambda_t.$ Such randomised
  strategies will not be used in this paper though. 

  Due to its flexibility, we envision manifold applications of the defective
  dynamics in our research in progress on polygonal fields. In the first place
  we are planning to use it to gain further knowledge on the higher order
  correlation structure of consistent polygonal fields. To illustrate the
  techniques we are currently developing we shall employ the defective 
  dynamics to establish exponential decay of correlations in the particular
  case of rectangular fields in Section \ref{CORRDEC} below. Another
  important field in which we anticipate the use of this dynamics is digital
  image segmentation \--- the defective dynamics will provide a rich and flexible
  class of new Monte-Carlo moves for polygonal field-based Bayesian segmentation
  algorithms we developed in joint work with Kluszczy\'nski and van Lieshout
  \cite{KLS05,KLS07,LS07,SL08}.
 
\section{Exponential decay of dependencies for consistent rectangular fields}\label{CORRDEC}
 The purpose of this section is to make use of the defective disagreement
 loop dynamics developed in Subsection \ref{DDLD} above in order to 
 establish exponential decay of higher order correlations for the consistent
 rectangular field ${\cal A}^{\cal M}.$ To simplify the presentation we assume
 throughout this section with no loss of generality that the admissible
 directions $l, l'$ of the field are parallel to the coordinate axes.
 Assume also that a collection 
 $(l_i,x_i)_{i=1}^k,\; x_i \in l_i,$ in general position as in Section
 \ref{CORCR} is given such that $x_1$ is its topmost point
 with the highest $y$ coordinate. Further, let ${\cal R}[x_2,\ldots,x_k]$ be
 the smallest rectangle with sides parallel to coordinate axes containing
 all remaining points $x_2,\ldots,x_k$ of the considered collection.
 We know from Lemma \ref{CORRHC} that if $l_1$ does not hit the convex
 hull of $x_2,\ldots,x_k$ then the correlation function
 $\sigma^{\cal M}[l_1,x_1;l_2,x_2;\ldots;l_k,x_k]$ factorises as
 $\sigma^{\cal M}[l_1,x_1] \sigma^{\cal M}[l_2,x_2;\ldots;l_k,x_k].$ 
 Here in Lemma \ref{EDec} we show that also when $l_1$ does hit
 $\conv(\{x_2,\ldots,x_k\})$ and hence it also hits the rectangle
 ${\cal R}[x_2,\ldots,x_k],$
 the correlation $\sigma^{\cal M}[l_1,x_1;l_2,x_2;\ldots;l_k,x_k]$ 
 approaches $\sigma^{\cal M}[l_1,x_1] \sigma^{\cal M}[l_2,x_2;\ldots;l_k,x_k]$
 exponentially fast in the distance measure between 
 $x_1$ and ${\cal R}[x_2,\ldots,x_k]$ along $l_1.$  
 \begin{lemma}\label{EDec}
  With notation as above, if $l_1$ hits ${\cal R}[x_2,\ldots,x_k]$ then we have
  $$ \sigma^{\cal M}[l_1,x_1;l_2,x_2;\ldots;l_k,x_k] 
    = \sigma^{\cal M}[l_1,x_1] \sigma^{\cal M}[l_2,x_2;\ldots;l_k,x_k] $$
  $$ \left(1+O\left(\e^{-{\cal M}([[l_1[x_1 \leftrightarrow {\cal R}[x_2,\ldots,x_k]]]])}
     \right)\right), $$
  where $l_1[x_1 \leftrightarrow {\cal R}[x_2,\ldots,x_k]]$ is the segment 
  of $l_1$ between $x_1$ and ${\cal R}[x_2,\ldots,x_k].$
 \end{lemma}

 \paragraph{Proof}
  Our proof relies on estimating the conditional probability of the event
  $$ {\cal E}[dl_1,x_1] := \{ \exists_{e \in \Edges({\cal A}^{\cal M})} \; 
      x_1 \in e,\; l[e] \in dl_1 \} $$
  given the event 
  $$ {\cal E}[d l_2,x_2;\ldots;d l_k,x_k] := 
     \{ \forall_{i=2}^k \exists_{e \in \Edges({\cal A}^{\cal M})} \;
        x_i \in e,\; l[e] \in dl_i \}. $$
  We shall show that 
  \begin{equation}\label{DOPOK}
   {\Bbb P}\left({\cal E}[dl_1,x_1] \left| {\cal E}[dl_2,x_2;\ldots;dl_k,x_k] \right. \right) =
     {\cal M}(dl_1)
     \left(1+O\left(\e^{-{\cal M}([[l_1[x_1 \leftrightarrow {\cal R}[x_2,\ldots,x_k]]]])}
     \right)\right), 
  \end{equation}
  which will complete the proof in view of (\ref{CORR}) [recall that $x_i \in l_i,\;
  i=1,\ldots,k\,$]
  and Theorem \ref{Corr12}.
  To proceed, denote by $\tilde{\cal A}^{\cal M}$ a rectangular field coinciding
  in law with ${\cal A}^{\cal M}$ conditioned on the event 
  ${\cal E}[d l_2,x_2;\ldots;d l_k,x_k].$
  Note that since $l_1$ crosses ${\cal R}[x_2,\ldots,x_k],$ it has to be a
  vertical line. Construct the following directional decision rule
  $\phi(\cdot,\cdot)$ for the defective dynamics:
  \begin{itemize}
   \item Whenever turning into a horisontal direction, move away from
         $l_1,$ that is to say rightwards if the turning point lies to
         the right from $l_1$ and leftwards if the turning point lies
         to the left of $l_1.$ When the turning point lies on $l_1,$
         move rightwards. 
   \item Whenever turning into a vertical direction, move upwards.
  \end{itemize}
  Moreover, let $D := {\cal R}[x_1,x_2,\ldots,x_k],$ that is to say $D$ is
  the smallest rectangle with sides parallel to coordinate axes, spanned by
  $x_1,\ldots,x_k.$ Furthermore, set $A_L$ to be the family of all vertical lines
  in $[[D]]$ and for each $l \in A_L = [[D]]$ let ${\Bbb A}^*(l)$ be the topmost
  point of $l \cap D$ so that in particular ${\Bbb A}^*(l_1) = x_1.$  
  With this notation, using Theorem \ref{DefAR1} including the reversibility
  statement there we see that the law of the conditional field $\tilde{\cal A}^{\cal M}_D$
  is invariant and reversible for the following conditional version of the
  defective disagreement loop dynamics:
  \begin{itemize}
   \item With intensity ${\cal M}(dl)ds$ for $l \in A_L$
    set $\gamma_{s+ds} := \gamma_s \oplus (l,{\Bbb A}^*(l)).$
    Should the proposed update fail or result in a configuration
    violating ${\cal E}[dl_2,x_2;\ldots;dl_k,x_k],$ 
    keep $\gamma_{s+ds} = \gamma_s.$
  \item For each line $l \in A_L$ extending an edge $e$ of 
   $\gamma_s$ such that ${\Bbb A}^*(l) \in e,$ with intensity $ds$
   set $\gamma_{s+ds} := \gamma_s \ominus (l,{\Bbb A}^*(l)).$
   Should the proposed update fail or result in a configuration
   violating ${\cal E}[dl_2,x_2;\ldots;dl_k,x_k],$ 
   keep $\gamma_{s+ds} = \gamma_s.$       
 \end{itemize}    
 Note that the only difference compared to the {\bf DefDL} dynamics 
 in Subsection \ref{DDLD} is that we discard updates resulting in
 violation of ${\cal E}[dl_2,x_2;\ldots;dl_k,x_k].$ Observe also
 that
 \begin{itemize}
  \item Under the so-defined directional decision rules the disagreement
        loops arising in the course of the dynamics are always
        single-branched (the second branch is chopped
        off by the boundary of $D$) and contain no cycles so that the
        updates in the above dynamics never fail (but note that they
        can nevertheless be discarded whenever violating 
        ${\cal E}[dl_2,x_2;\ldots;dl_k,x_k]$).
  \item The only update of the dynamics which can lead to a transition from a 
        configuration where ${\cal E}[dl_1,x_1]$ holds to a configuration where
        it does not hold is $\ominus (l_1,x_1).$ Indeed, in view of the directional
        decision rules, even if some disagreement loop annihilates a subsegment of
        $l_1,$ this subsegment cannot reach $x_1.$ 
  \item Likewise, the only update of the dynamics which can lead to a transition from
        a configuration where ${\cal E}[dl_1,x_1]$ does not hold to a configuration
        where it holds is $\oplus (l_1,x_1).$   
  \item For the update $\oplus (l_1,x_1)$ the only possibility to cause
        violation of ${\cal E}[dl_2,x_2;\ldots;dl_k,x_k],$ and thus to
        be discarded, is that the initial downwards creation phase growth
        of the disagreement loop from $x_1$ continue with no directional
        updates along $l_1$ until the loop reaches ${\cal R}[x_2,\ldots,x_k].$ Indeed,
        if a directional update occured prior to that, the disagreement
        loop would only unfold upwards and rightwards thereafter, never
        reaching ${\cal R}[x_2,\ldots,x_k]$ and consequently not violating
        ${\cal E}[dl_2,x_2;\ldots;dl_k,x_k].$ However, the probability
        that the initial segment of the disagreement loop for $\oplus (x_1,l_1)$
        extends all the way down to ${\cal R}[x_2,\ldots,x_k]$ without directional
        updates is $\exp(-{\cal M}([[l_1[x_1 \leftrightarrow
        {\cal R}[x_2,\ldots,x_k]]]]))$ by the definition of the dynamics.
        Clearly, even if this happens, the violation of ${\cal E}[dl_2,x_2;\ldots;dl_k,x_k]$
        is not guaranteed, but possible. 
  \item Likewise, the update $\ominus (l_1,x_1),$ only possible on the
        event ${\cal E}[dl_1,x_1],$ may only cause violation of
        ${\cal E}[dl_2,x_2;\ldots;dl_k,x_k]$ if the current configuration
        $\gamma_s$ contains a segment along $l_1$ joining $x_1$ to 
        ${\cal R}[x_2,\ldots,x_k].$ On the event 
        ${\cal E}[dl_1,x_1] \cap {\cal E}[dl_2,x_2;\ldots;dl_k,x_k]$
        this only happens with probability 
        $O(\exp(-{\cal M}([[l_1[x_1 \leftrightarrow {\cal R}[x_2,\ldots,x_k]]]])))$
        though, in analogy to the previous case.
 \end{itemize}
 In view of these observations and by the definition of the dynamics 
 we see that in its course:
 \begin{itemize}
  \item On $\neg {\cal E}[dl_1,x_1],$ transitions leading to ${\cal E}[dl_1,x_1]$
        occur in s-time with intensity
        $$ {\cal M}(dl_1) ds
           \left(1+O\left(\e^{-{\cal M}([[l_1[x_1 \leftrightarrow {\cal R}[x_2,\ldots,x_k]]]])}
           \right)\right).$$
  \item On ${\cal E}[dl_1,x_1],$ transitions leading to $\neg {\cal E}[dl_1,x_1]$
        occur in s-time with intensity
        $$ ds \left(1+O\left(\e^{-{\cal M}([[l_1[x_1 \leftrightarrow {\cal R}[x_2,\ldots,x_k]]]])}
           \right)\right). $$
 \end{itemize}
 Consequently, using the detailed balance we see that the probability of 
 ${\cal E}[dl_1,x_1]$ for $\tilde{\cal A}^{\cal M},$ or equivalently the 
 conditional probability of ${\cal E}[dl_1,x_1]$ given
 ${\cal E}[dl_2,x_2;\ldots;dl_k,x_k],$ is 
 $$ {\cal M}(dl_1) \left(1+O\left(\e^{-{\cal M}([[l_1[x_1 \leftrightarrow {\cal R}[x_2,\ldots,x_k]]]])}
     \right)\right) $$
 as required for (\ref{DOPOK}). This completes
 the proof of the lemma. $\Box$

 A repetitive use of Lemma \ref{EDec} readily yields the following approximate factorisation theorem.
 \begin{theorem}\label{EDecT}
  Assume the collection $(l_i,x_i)_{i=1}^k$ is such that $x_i \not \in {\cal R}[x_{i+1},\ldots,x_k]$
  for all $i$ and let $\delta_i := {\cal M}[[l_i[x_i \leftrightarrow {\cal R}[x_{i+1},\ldots,x_k]]]]$
  if $l_i$ hits ${\cal R}[x_{i+1},\ldots,x_k]$ and $\delta_i := +\infty$ otherwise.
  Then
  $$ \sigma^{\cal M}[dl_1,x_1;\ldots;dl_k,x_k] = \left(\prod_{i=1}^k {\cal M}(dl_i)\right)
     \prod_{i=1}^k (1+O(\exp(-\delta_i))). $$
 \end{theorem}

\section{Contour birth and death dynamics}\label{CBDD}
 The purpose of the present section is to describe an alternative dynamics
 which can be used to simulate the empty-boundary field 
 ${\cal A}^{{\cal M};\beta}_{D|\emptyset}$ for $\beta \geq 2.$ 
 A particular feature of this dynamics is that it yields a perfect sampler, i.e.
 the sample it outputs comes exactly from the target distribution and not just
 its approximation. Moreover, for $\beta$ large enough it extends to the whole
 plane thus allowing for a direct graphical construction (perfect simulation)
 of the whole plane thermodynamic limit ${\cal A}^{{\cal M};\beta}.$ The idea
 underlying the graphical representation comes from Fern\'andez, Ferrari \& Garcia
 \cite{FFG1,FFG2,FFG3} and in the context of homogeneous polygonal fields
 it has been developed in Schreiber \cite{SC05,SC06}. Here we present an
 extension of the construction to the general non-homogeneous setting. 
 Our presentation splits into several steps.
 \subsection{Free contour measure}
 Choose an open bounded set $D$ with piecewise smooth boundary and
 consider the space ${\cal C}_D$ consisting of all closed polygonal contours in $D$
 which do not touch the boundary $\partial D.$ For a given finite collection
 $(l) := ( l_1,...,l_n )$ of straight lines intersecting $D$ denote by ${\cal C}_D(l)$
 the family of those polygonal contours in ${\cal C}_D$ which belong to 
 $\Gamma_{D|\emptyset}(l) := \Gamma_D(l) \cap \Gamma_{D|\emptyset}.$
 We define the so-called free contour measure $\Theta_D$ on ${\cal C}_D$ by
 putting for $C \subseteq {\cal C}_D$ measurable, say with respect to the
 Borel $\sigma$-field generated by the Hausdorff distance topology,
 \begin{equation}\label{WOLNEKONTURY}
   \Theta_D(C) = \int_{{\rm Fin}(L[D])} \sum_{\theta \in C \cap {\cal C}_D(l)}
    \exp(-L^{\cal M}(\theta)) d{\cal M}^*((l))
 \end{equation}
 with ${\rm Fin}(L[D])$ standing for the for the family of finite
 collections of lines intersecting $D$ and where ${\cal M}^*$ is the measure
 on ${\rm Fin}(L[D])$ given by 
 $d{\cal M}^*((l_1,...,l_n)) := \prod_{i=1}^n {\cal M}(dl_i).$
 Thus, in other words,
 $$ \Theta_D(d\theta) = \exp(-L^{\cal M}(\theta)) \prod_{e\in\Edges(\theta)} {\cal M}(dl[e]). $$
 For general $\beta \in {\Bbb R}$ we consider the exponential modification
 $\Theta^{[\beta]}_D$ of the free measure $\Theta_D,$
 \begin{equation}\label{THETAB}
  \Theta^{[\beta]}_D(d\theta) := \exp(-[\beta-1] L^{\cal M}(\theta))
  \Theta_D(d\theta) = \exp(-\beta L^{\cal M}(\theta)) \prod_{e\in\Edges(\theta)} {\cal M}(dl[e]).
 \end{equation}
 It is easily seen that the total mass $\Theta^{[\beta]}_D({\cal C}_D)$ is always
 finite for bounded $D$ for all $\beta \in {\Bbb R}.$
 Indeed, taking into account that $L^{\cal M}(\theta) \leq
 \card(\Edges(\theta)) {\cal M}([[D]])$ we have in view of (\ref{THETAB})
 \begin{equation}\label{THSK}
   \Theta^{[\beta]}_D({\cal C}_D) \leq \sum_{k=0}^{\infty} \frac{1}{k!}
    [{\cal M}([[D]])]^k \exp(|\beta|k {\cal M}([[D]])) = \exp({\cal M}([[D]])
    \exp(|\beta|{\cal M}([[D]]))) < \infty. 
 \end{equation}
 Note that for all $D$ the free contour measures $\Theta_D$ as defined
 in (\ref{WOLNEKONTURY}) arise as the respective restrictions to ${\cal C}_D$
 of the same measure $\Theta$ on ${\cal C} := \bigcup_{n=1}^{\infty}
  {\cal C}_{(-n,n)^2},$ in the sequel referred to as the infinite volume
 free contour measure. Indeed, this follows easily by the observation
 that $\Theta_{D_1}$ restricted to ${\cal C}_{D_2}$ coincides with
 $\Theta_{D_2}$ for $D_2 \subseteq D_1.$ In the same way we construct
 the infinite-volume exponentially modified measures 
 $\Theta^{[\beta]},\; \beta \in {\Bbb R}.$ It turns out that for
 $\beta \geq 2$ a natural algorithmic random walk representation of 
 $\Theta^{[\beta]}$ can be given. This is done as follows. For each
 contour $\theta \in {\cal C}$ let $\iota[\theta]$ be its leftmost
 vertex, that is to say its extreme left point minimising the first 
 coordinate, with possible ties broken in an arbitrary measurable 
 way. Next, for $x \in {\Bbb R}^2$ denote by ${\cal C}_x$ the collection of
 contours $\theta \in {\cal C}$ such that $\iota[\theta] = x.$ For
 $\beta \geq 2$ let $\Theta_x^{[\beta]}$ be the subprobability
 measure determined by the following construction of a ${\cal C}_x$-valued
 $\Theta_x^{[\beta]}$-distributed random element $\theta$
\begin{itemize}
  \item Simulate a continuous-time random walk $(Z_t)_{t \geq 0}$ gouverned by the
        following dynamics
       \begin{itemize}
        \item choose two random lines $l_1,l_2$ meeting at $x$ according
              to the distribution $$\frac{d[{\cal M} \times {\cal M}]
              ((l_1,l_2),\; l_1 \cap l_2 \in dx)}{2 \lnaw {\cal M}\pnaw(dx)},$$
              as in (\ref{THETAWPRO}), 
        \item put $Z_0 := 0$ and choose one of two possible initial directions along $l_1$
              with equal probabilities $1/2,$ 
        \item between direction update events specified below move in a constant direction with 
              speed $1,$
        \item while moving along the segment $Z_{[t,t+dt]} = \overline{Z_t Z_{t+dt}},$
              update the movement direction and start moving along a line $l \in [[Z_{[t,t+dt]}]]$
              with probability $2 {\cal M}(dl)$ and with equiprobable choices between the
              two possible directions along $l,$  
       \end{itemize}
  \item Consider a killed modification $\tilde{Z}^{[\beta-2]}_t$ of $Z_t$ by
        killing $Z_t$
        \begin{itemize}
         \item with intensity $[\beta-2] {\cal M}([[Z_{[t,t+dt]}]]),$
         \item whenever $Z_t$ hits its past trajectory,
        \end{itemize} 
  \item Draw an infinite {\it loop-closing} half-line $l^*$ along $l_2$ beginning at $x$
        with equiprobable choices of two possible directions,  
  \item If the random walk $\tilde{Z}^{[\beta-2]}_t$ hits the loop closing line $l^*$
        before getting killed, and the self-avoiding contour 
        $\theta_* := \theta_*[\tilde{Z}^{[\beta-2]};l^*]$ created by
        $l^*$ and the trajectory of $\tilde{Z}^{[\beta-2]}_t$ up to the moment
        of hitting $l^*$ satisfies $\iota[\theta_*] = x,$ then
        \begin{itemize}
         \item with probability $\exp(-\beta{\cal M}([[e^*]]))$ output
               $\theta := \theta_*,$ where $e^*$ stands for the
               segment of the loop closing line $l^*$ from $x$ to its
               intersection point with $\tilde{Z}^{[\beta-2]}_t,$
         \item otherwise output $\theta := \emptyset.$
        \end{itemize}
        In all remaining cases set $\theta := \emptyset.$
 \end{itemize}        
 Note that outputting $\emptyset$ means failing in the above algorithm,
 and hence the resulting distribution is a sub-probability rather than
 probability measure. 

 The following lemma provides a constructive representation for
 $\Theta^{[\beta]},\; \beta \geq 2,$ and is close in spirit to
 Lemma 5.1 in \cite{SC06} and Lemma 1 in \cite{LS07}. 
\begin{lemma}\label{MPPrepr}
 For $\beta \geq 2$ we have
 $$ \Theta^{[\beta]} = 4 \int_{{\Bbb R}^2} \Theta_x^{[\beta]} \lnaw{\cal M}\pnaw(dx). $$ 
\end{lemma}

\paragraph{Proof}
 The directed nature of the random walk trajectories as constructed above requires
 considering for each contour $\theta$ two oriented instances $\theta^{\rightarrow}$
 (clockwise) and $\theta^{\leftarrow}$ (anti-clockwise). By the construction 
 of $\Theta_x^{[\beta]}$ the assertion of the lemma will follow as soon as we show
 that for each $x \in {\Bbb R}^2$ and $\theta \in {\cal C}_x$ we have 
 \begin{equation}\label{RWR}
  8 \lnaw {\cal M} \pnaw (dx) \e^{-\beta L^{\cal M}(e^*)} 
  {\Bbb P}\left( \tilde{Z}^{[\beta-2]}_t \mbox{ reaches $l^*$ and }
  \theta_*[\tilde{Z}^{[\beta-2]};l^*] \in d \theta^{\rightarrow} \right) = 
  \Theta^{[\beta]}(d\theta),
\end{equation} 
where $e^*$ stands for the last segment of $\theta^{\rightarrow}$ counting from $x$ as 
the initial vertex, which is to coincide with the segment of the loop-closing line $l^*$ 
joining to $x$ its intersection point with $\tilde{Z}^{[\beta-2]}_t;$ whereas 
$\theta_*[\tilde{Z}^{[\beta-2]};l^*]$ is the self-avoiding contour created by
$l^*$ and the trajectory of $\tilde{Z}^{[\beta-2]}_t$ up to the moment of hitting $l^*,$
as denoted in the construction of $\Theta_x^{[\beta]}$ above. Indeed, the same 
relation holds then for $\theta^{\leftarrow},$ so adding versions of (\ref{RWR}) for 
$\theta^{\rightarrow}$ and $\theta^{\leftarrow},$ which amounts to taking into account 
two possible directions in which the random walk can move along $\theta,$ will yield 
$2\Theta^{[\beta]}(d\theta)$ on the right hand side, whence the assertion of the
lemma will follow. 

To establish (\ref{RWR}), we observe that the probability element 
$$
{\Bbb P}\left( \tilde{Z}^{[\beta-2]}_t \mbox{ reaches $l^*$ and }
 \theta_*[\tilde{Z}^{[\beta-2]};l^*] \in  d \theta^{\rightarrow} \right)
$$
is exactly 
\begin{equation}\label{RWR2}
   \frac{1}{8 \lnaw {\cal M} \pnaw (dx)}
   \exp\left(-\beta L^{\cal M}(\theta \setminus e^*)\right) \prod_{i=1}^k {\cal M}(dl[e_k]),
\end{equation}
where $e_1,\ldots,e_k$ are all segments of $\theta$ including $e^* = e_k,$
while $l[e_i]$ stands for the straight line determined by $e_i.$ Indeed,
 \begin{itemize}
  \item the prefactor $[8 \lnaw {\cal M} \pnaw (dx)]^{-1} = 
        [4 [{\cal M}\times {\cal M}](\{(l_1,l_2)\;|\; l_1 \cap l_2 \in dx \})]^{-1}$
        times the product ${\cal M}(dl[e_1]) {\cal M}(dl[e_k])$ comes from the
        choice of the lines $l_1$ and
        $l_2$ containing respectively the initial segment of $\theta^{\rightarrow}$
        (counting from $x$) and the half-line $l^*,$ as well as from the
        choice between two equiprobable
        directions on each of these lines,
  \item for each of the remaining segments $e_i,\; i=2,\ldots,k-1,$ the factor
        ${\cal M}(dl[e_i])$ comes from the directional update of the random 
        walk $Z_t,$ with $2 {\cal M}(dl[e_i])$ due to the choice of the line
        and $1/2$ due to the choice between two equiprobable directions along
        this line,  
  \item finally, $\exp(-[\beta-2] L^{\cal M}(\theta \setminus e^*))$ comes
        from killing in the course of the random walk whereas extra
        $\exp(-2 L^{\cal M}(\theta \setminus e^*))$ is due to the absence
        of directional updates along the segments of the walk, which
        yields $\exp(-\beta L^{\cal M}(\theta \setminus e^*))$ when put
        together.  
  \end{itemize}
 Recalling (\ref{THETAB}) we see that the expression in (\ref{RWR2}) coincides with
 $$\frac{1}{8 \lnaw {\cal M} \pnaw (dx)} \exp(\beta L^{\cal M}(e^*))
  \Theta^{[\beta]}(d\theta),$$ which yields the required relation (\ref{RWR}).
 The proof is complete. $\Box$

 To proceed, assume that the activity measure ${\cal M}$ admits a homogeneous upper
 and lower bound, that is to say
 \begin{equation}\label{HOMOGENBD}
  C_- \mu \leq {\cal M} \leq C_+ \mu 
 \end{equation}
 for some $0 < C_- \leq C_+ < \infty$ and with $\mu$ standing for the standard
 isometry-invariant Haar-Lebesgue measure on the space $[[{\Bbb R}^2]]$ of straight lines
 in ${\Bbb R}^2.$ Recall that one possible construction of $\mu$ goes by identifying a straight
 line $l$ with the pair $(\phi,\rho) \in [0,\pi) \times {\Bbb R},$
 where $(\rho \sin(\phi), \rho \cos(\phi))$ is the vector orthogonal to $l,$ and joining it to the
 origin, and then by endowing the parameter space $[0,\pi) \times {\Bbb R}$ with the usual Lebesgue
 measure. Then, for each linear segment $e$ we have by standard integral geometry 
 \begin{equation}\label{HH1}
  2 C_- \ell(e) \leq {\cal M}([[e]]) \leq 2 C_+ \ell(e)
 \end{equation}
 because $\mu([[e]]) = 2 \ell(e).$
 Likewise, since $\lnaw \mu \pnaw (dx) = \pi dx$ as easily checked by straightforward
 integration,  
 \begin{equation}\label{HH2}
  C_- \pi dx \leq \lnaw {\cal M} \pnaw (dx) \leq C_+ \pi dx.
 \end{equation}
 Consequently, we conclude from Lemma \ref{MPPrepr} that under
 $\Theta^{[\beta]},\; \beta \geq 2,$ the contour size exhibits
 exponentially decaying tails, which is a non-homogeneous counterpart
 of Lemma 1 in \cite{SC05}.  
 \begin{lemma}\label{NICHOLLS}
  Under the assumption (\ref{HOMOGENBD}), for $\beta \geq 2$ we have
  $$
    \Theta^{[\beta]}(\{ \theta;\; dx \cap \Vertices(\theta) \neq \emptyset,\; \lgth(\theta) > R \})
    \leq 
    4 C_+ \pi \exp(- 2 C_- [\beta-2] R) dx,
  $$
  where the event $\{ dx \cap \Vertices(\theta) \neq \emptyset \}$ is to be understood
  that a vertex of $\theta$ falls into $dx.$ 
  Moreover, there exists a constant $\varepsilon > 0$ such that, for $\beta \geq 2,$
  $$
   \Theta^{[\beta]}(\{ \theta;\; {\bf 0} \in \innt \theta,\;
   \lgth(\theta) > R \}) \leq \exp(-2 C_- [\beta-2+\varepsilon] R + o(R)),
  $$
  with $\innt \theta$ standing for the region enclosed by the contour $\theta.$
 \end{lemma}
 \paragraph{Proof} The first assertion follows by the construction of the random walk $Z_t$ 
  where, in view of (\ref{HH1}), the killing intensity is at least $2 C_-$ times the length
  element covered under the present assumptions, whereas the extreme left vertices of
  contours have their intensity bounded by $4 \lnaw {\cal M} \pnaw(dx) \leq 4 C_+ \pi$
  in view of the integral formula in Theorem \ref{MPPrepr} and of (\ref{HH2}). To get the
  second assertion observe in addition that during each unit time of its evolution the
  random walk $Z_t$ has some uniformly non-zero chance of hitting its past trajectory,
  see also the proof of Lemma 1 in \cite{SC05}. $\Box$  
    
 \subsection{Polymer representation}
  To proceed we let ${\cal P}_{\Theta_D^{[\beta]}}$ be the Poisson point process in
  ${\cal C}_D$ with intensity measure $\Theta^{[\beta]}_D.$ It follows then directly
  by (\ref{THETAB}) in view of (\ref{GIBBS1}) and (\ref{GIBBS2}) that, for all $\beta \in {\Bbb R},$ 
  ${\cal A}^{{\cal M};\beta}_{D|\emptyset}$
  coincides in distribution with the union of contours in ${\cal P}_{\Theta_D^{[\beta]}}$
  conditioned on the event that they are disjoint so that
  \begin{equation}\label{WARUNKOWKT}
   {\cal L}\left( {\cal A}^{{\cal M};\beta}_{D|\emptyset} \right) =
   {\cal L}\left( \bigcup_{\theta \in {\cal P}_{\Theta_D^{[\beta]}}} \theta \; \left|
   \; \forall_{\theta, \theta' \in {\cal P}_{\Theta_D^{[\beta]}}} \theta \neq \theta'
   \Rightarrow \theta \cap \theta' = \emptyset \right. \right),
 \end{equation}
 where the conditioning makes sense because $\Theta_D^{[\beta]}({\cal C}_D)$ is
 finite as shown in (\ref{THSK}) above, see also \cite{SC05}, Section 2.2.
 In particular, in analogy to Subsection 2.2 and Theorem 2 ibidem,
 the law of ${\cal A}^{{\cal M};\beta}_{D|\emptyset}$ is invariant and reversible with
 respect to the following contour birth and death dynamics $(\gamma^D_s)$ on $\Gamma_D.$
 \begin{description}
  \item{${\bf (C:birth[\beta])}$} With intensity $\Theta^{[\beta]}_D(d\theta) ds$ do
   \begin{itemize}
    \item Choose a new contour $\theta,$
    \item If $\theta \cap \gamma^D_s = \emptyset,$ accept $\theta$ and
          set $\gamma^D_{s+ds} := \gamma^D_s \cup \theta,$
    \item Otherwise reject $\theta$ and keep $\gamma^D_{s+ds} := \gamma^D_s,$
   \end{itemize}
  \item{${\bf (C:death[\beta])}$}
      With intensity $1 \cdot ds$ for each contour $\theta \in \gamma^D_s$
      remove $\theta$ from $\gamma^D_s$ setting $\gamma^D_s := \gamma^D_s \setminus \theta.$
 \end{description}
 Moreover, ${\cal L}({\cal A}^{{\cal M};\beta}_D)$ is the unique invariant distribution
 of the above dynamics, see Theorem 2 in \cite{SC05}. Clearly, both the representation
 (\ref{WARUNKOWKT}) and the above dynamics can be regarded as a kind of polymer representation
 for ${\cal A}^{{\cal M};\beta}_D.$

 Unlike the results for free contour measures, the polymer representation is valid
 for all values of $\beta.$ In this context it is helpful to note that, by (\ref{WARUNKOWKT}),
 the probability ${\Bbb P}({\cal A}^{{\cal M};\beta}_{D|\emptyset} = \emptyset)$ is not smaller than 
 $\exp(-\Theta^{[\beta]}_D({\cal C}_D)).$ On the other hand, comparing with (\ref{GIBBS1})
 shows that ${\Bbb P}({\cal A}^{{\cal M};\beta}_{D|\emptyset} = \emptyset)$ is 
 $[{\cal Z}^{{\cal M};\beta}_{D|\emptyset}]^{-1}$ where, recall, 
 ${\cal Z}^{{\cal M};\beta}_{D|\emptyset}$ is the partition function for 
 ${\cal A}^{\cal M}_{D|\emptyset}.$ Thus, by (\ref{THSK}), 
 \begin{equation}\label{ZZ1}
  {\cal Z}^{{\cal M};\beta}_{D|\emptyset} \leq \exp(\Theta^{[\beta]}({\cal C}_D)) < \infty.
 \end{equation}
 Likewise,
 \begin{equation}\label{ZZ2}
  {\cal Z}^{{\cal M};\beta}_D < \infty 
 \end{equation}
 which can be proven along the same lines, see also the proof of Corollary 2 in \cite{SC05}.

 \subsection{Graphical representation}
 These observations
 place us within the framework of the general contour birth and death graphical
 construction as developed by Fern\'andez, Ferrari \& Garcia
 \cite{FFG1,FFG2,FFG3} and as sketched below, see ibidem and \cite{SC05} for
 further details. Choose $\beta$
 large enough, to be specified below.  Define ${\cal F}({\cal C})$ to be the space of
 countable and locally finite collections of contours from ${\cal C},$ with
 the local finiteness requirement meaning that at most a finite number of
 contours can hit a bounded subset of ${\Bbb R}^2.$ On the time-space
 ${\Bbb R} \times {\cal F}({\cal C})$ we construct the time-stationary free contour
 birth and death process $(\varrho_s)_{s \in {\Bbb R}}$ with the birth intensity
 measure given by $\Theta^{[\beta]}$ and with the death intensity $1.$
 Note that {\it free} means here that every new-born contour is accepted
 regardless of whether it hits the union of already existing contours or not,
 moreover we admit negative time here, letting $s$ range through ${\Bbb R}$
 rather than just ${\Bbb R}_+.$
 Observe also that we need the birth measure $\Theta^{[\beta]}$ to be finite on
 the sets $\{ \theta \in {\cal C} \;|\; \theta \cap A \neq \emptyset \}$
 for all bounded Borel $A \subseteq {\Bbb R}^2$ in order to have the process
 $(\varrho_s)_{s \in {\Bbb R}}$ well defined on ${\Bbb R} \times {\cal F}({\cal C}).$
 By Lemma \ref{NICHOLLS} this is ensured whenever $\beta \geq 2.$
 To proceed, for the free process $(\varrho_s)_{s \in {\Bbb R}}$ we perform the
 following {\it trimming} procedure. We place a directed connection from each
 time-space instance of a contour showing up in $(\varrho_s)_{s \in {\Bbb R}}$
 and denoted by $\theta \times [s_0,s_1),$ with $\theta$ standing for the contour
 and $[s_0,s_1)$ for its lifespan, to all time-space contour instances
 $\theta' \times [s'_0,s'_1)$ with $\theta' \cap \theta \neq \emptyset,
 s'_0 \leq s_0$ and $s'_1 > s_0.$ In other words, we connect $\theta \times
 [s_0,s_1)$ to its {\it ancestors} understood as those contour instances which may have
 affected the acceptance status of $\theta \times [s_0,s_1)$ in the {\it constrained}
 contour birth and death dynamics {\bf (C)} as discussed above.
 These directed connections give rise to directed ancestor chains of time-space
 contour instances, following Fernandez, Ferrari \& Garcia \cite{FFG3} the union
 of all ancestor chains stemming from a given contour instance $\theta^* =
 \theta \times [s_0,s_1),$ including the instance itself,
 is referred to as its {\it clan of ancestors}
 and is denoted by $\An(\theta^*).$ More generally, for a bounded region
 $U$ in the plane we write $\An_s(U)$ for the union of ancestor clans of
 all contour instances $\theta \times [s_0,s_1)$ with $\theta \cap U \neq
 \emptyset$ and $s \in [s_0,s_1).$ Lemma \ref{NICHOLLS} allows us to apply
 the technique of domination by sub-critical branching processes, developed
 in \cite{FFG1,FFG2,FFG3}, in order to conclude
 that there exists $\beta_g$ such that for each $\beta > \beta_g$ there
 exists $c := c(\beta) > 0$ such that
 \begin{equation}\label{ZANIKGR}
  {\Bbb P}(\diam \An_s({\Bbb B}_2(x,1)) > R) \leq \exp(- c R),\; s \in {\Bbb R},\; x \in {\Bbb R}^2,
 \end{equation}
 with ${\Bbb B}_2(x,1)$ standing for the ball of radius $1$ in ${\Bbb R}^2$
 centred at $x.$ In the sequel we shall always assume that $\beta > \beta_g,$
 that is to say that $\beta$ is in the {\it validity} region
 of the graphical construction. We see that for $\beta > \beta_g$ all the
 ancestor clans are a.s. finite and we can uniquely determine the acceptance
 status of all their members: contour instances with no ancestors are a.s.
 accepted, which automatically and uniquely determines the acceptance status
 of all the remaining members of the clan by recursive application of
 the inter-contour exclusion rule.
 In this case, discarding the unaccepted contour instances leaves us with a
 time-space representation of a time-stationary evolution $(\gamma_s)_{s \in {\Bbb R}}$
 on ${\cal F}({\cal C}),$ which is easily checked to evolve according to the
 whole-plane version of the dynamics {\bf (C)} above. In full analogy with Section 4
 and Theorem 4 of \cite{SC05} we see that for all $s \in {\Bbb R}$ the polygonal field
 $\gamma_s$ coincides in distribution with the thermodynamic limit (see Section
 3 ibidem) for ${\cal A}^{{\cal M};\beta}$ without infinite contours, which is unique
 (see Corollary 4 ibidem). For definiteness we put ${\cal A}^{{\cal M};\beta} := \gamma_0.$
 It should be also observed that for each $s \in {\Bbb R}$ the free field $\varrho_s$
 coincides in distribution with the Poisson contour process ${\cal P}_{\Theta^{[\beta]}}.$
 Since almost surely we have $\gamma_s \subseteq \varrho_s,$ we get the stochastic
 domination of the contour ensemble ${\cal A}^{{\cal M};\beta}$ by ${\cal P}_{\Theta^{[\beta]}}.$
 Moreover, using the exponential decay relation (\ref{ZANIKGR}), in full analogy with
 Theorem 4 (iv,v) we readily establish the following $\beta$-mixing statement for
 ${\cal A}^{{\cal M};\beta}.$
 \begin{corollary}\label{BetaMix}
  For $\beta > \beta_g$ there exists $c > 0$ such that for each $A,B \subseteq {\Bbb R}^2$
  the law of the restriction ${\cal A}^{{\cal M};\beta}_{|A \cup B}$ of the thermodynamic
  limit ${\cal A}^{{\cal M};\beta}$ to $A \cup B$ differs in variational distance
  from the product of the laws of the respective restrictions ${\cal A}^{{\cal M};\beta}_{|A},\;
  {\cal A}^{{\cal M};\beta}_{|B}$ by at most $O({\rm cover}(A) {\rm cover}(B) \exp(-c \dist(A,B))),$
  where ${\rm cover}(X)$ stands for the minimum number of unit balls covering $X \subseteq {\Bbb R}^2$
  whereas $\dist(A,B) := \inf_{x \in A, y \in B} \dist(x,y).$
 \end{corollary} 
 Indeed, this is easily checked by noting that the considered variational distance
 is bounded above by the probability that either $\An_0(A)$ extends further than
 $\dist(A,B)/2$ away from $A$ or $\An_0(B)$ extends further than $\dist(A,B)/2$
 away from $B$ which happens with probability at most $O({\rm cover}(A) {\rm cover}(B)
 \exp(-c \dist(A,B)))$ as required.
 Furthermore, in analogy to the discussion following Corollary 4 in \cite{SC05}
 we can use the stochastic domination of ${\cal A}^{{\cal M};\beta}$ by 
 ${\cal P}_{\Theta^{[\beta]}}$ and apply Lemma  \ref{NICHOLLS} to conclude that
 \begin{corollary}\label{FinContNest}
  With probability $1$ there is only finite contour nesting in ${\cal A}^{{\cal M};\beta},$
  that is to say no infinite chains of nested contours are present.
 \end{corollary}
 Indeed, this comes from the fact that the expected number of contours surrounding a 
 given point is always finite by Lemma \ref{NICHOLLS}.
 In addition, for $\beta$ large enough this expected
 number can be made arbitrarily small, which clearly guarantees the presence of 
 long-range point-to-point correlations and spontaneous magnetisation with contours
 interpreted as separating phases of different signs, as follows by a standard
 Peierls-type argument.

 We also consider finite-volume versions of the above
 graphical construction, replacing the infinite-volume birth intensity measure
 $\Theta^{[\beta]}$ with its finite-volume counterparts $\Theta^{[\beta]}_D$
 for bounded and open $D$
 with piecewise smooth boundary. Clearly, the graphical construction yields
 then a version of the finite-volume contour birth and death evolution ${\bf (C)}.$
 For each $D$ denote by $(\gamma_s^D)_{s \in {\Bbb R}}$ the resulting finite-volume
 time-stationary process on the space ${\cal F}({\cal C}_D)$ of finite contour configurations
 in $D$ and write $(\varrho^D_s)_{s \in {\Bbb R}}$ for the corresponding free process.
 It follows as in Theorem 2 in \cite{SC05} that $\gamma_s^D$
 coincides in distribution with ${\cal A}^{{\cal M};\beta}_D$ for each
 $s \in {\Bbb R}.$ Likewise, $\varrho^D_s$ coincides in distribution with
 ${\cal P}_{\Theta_D^{[\beta]}}.$

 By representing the measures $\Theta^{[\beta]}_D$ as the corresponding restrictions
 of $\Theta^{[\beta]}$ we obtain a natural coupling of all the processes
 $\gamma^D_s, \varrho^D_s, \gamma_s$ and $\varrho_s$ on a common probability space.
 We shall also consider ${\cal A}_D^{{\cal M};\beta}$ coupled on the same probability
 space by putting ${\cal A}_D^{{\cal M};\beta} = \gamma^D_0.$ Under this coupling
 we have the following corollary stating exponential convergence of 
 ${\cal A}^{{\cal M};\beta}_D$ to ${\cal A}^{{\cal M};\beta}$ as $D \uparrow {\Bbb R}^2.$
 \begin{corollary}\label{ZBIEGD}
  Assume that $\beta$ is within the validity regime of the contour birth and
  death graphical construction. Then for $x \in D$ such that $B(x,1) \subseteq D$
  we have
  $$ {\Bbb P}\left( [{\cal A}^{{\cal M};\beta}]_{|{\Bbb B}_2(x,1)} \neq
                    [{\cal A}^{{\cal M};\beta}_D]_{|{\Bbb B}_2(x,1)} \right) \leq
     \exp(-c [\dist(x,\partial D)-1]), $$
  with $c$ as in (\ref{ZANIKGR}) and with $|$ standing for the restriction operator.
 \end{corollary}
  Indeed, this follows easily by (\ref{ZANIKGR}) upon noting that the event
  $\{[{\cal A}^{{\cal M};\beta}]_{|{\Bbb B}_2(x,1)} \neq
                    [{\cal A}^{{\cal M};\beta}_D]_{|{\Bbb B}_2(x,1)} \}$ only occurs
  when $\An_0({\Bbb B}_2(x,1)) \not\subseteq D.$  


\end{document}